\documentclass{amsproc}
\usepackage{graphicx,amssymb}

\newtheorem{theorem}{Theorem}[section]
\newtheorem{lemma}[theorem]{Lemma}
\theoremstyle{definition}
\newtheorem{definition}[theorem]{Definition}
\newtheorem{corollary}[theorem]{Corollary}

\theoremstyle{remark}
\newtheorem{remark}[theorem]{Remark}
\newtheorem{step}[theorem]{}

\graphicspath{{figs/}}

\def\p{\prime}

\def\K{\mathcal{K}}
\def\R{\mathbb{R}}
\def\Z{\mathbb{Z}}

\let\mgp=\marginpar \marginparwidth18mm \marginparsep1mm
\def\marginpar#1{\mgp{\raggedright\tiny #1}}

\let\lbl=\label
\def\label#1{\lbl{#1}\ifinner\else\marginpar{\ref{#1} #1}\ignorespaces\fi}

\makeatletter
\newbox\overbox
\def\fakeover#1{\setbox\overbox\hbox{$#1$}\hbox
            {$\overline{#1\hskip-\wd\overbox}$\hskip\wd\overbox}}

\def\overrightarrow#1{\mskip2mu\vbox{\m@th\ialign{##\crcr
\rightarrowfill\crcr

\noalign{\kern-.4pt               %compensate for thickness of rule
\kern-\fontdimen22\textfont2 %pull down arrow so axis is on baseline
\nointerlineskip}
${\mskip0mu\hfil\fakeover{#1}\hfil\mskip6mu}$\crcr}}\mskip-2mu}

\def\overleftrightarrow#1{\mskip-3mu\vbox{\m@th\ialign{##\crcr
\leftarrowfill\hskip-.6em\rightarrowfill\crcr
\noalign{\kern-.4pt               %compensate for thickness of rule
\kern-\fontdimen22\textfont2 %pull down arrow so axis is on baseline
\nointerlineskip}
${\mskip3mu\hfil\fakeover{#1}\hfil\mskip6mu}$\crcr}}\mskip-2mu}

\makeatother

\begin{document}

\title[Ropelengths of Knots are Almost Linear]{The Ropelengths of Knots Are Almost Linear in Terms of Their Crossing Numbers}


\author{Yuanan Diao}
\address{Department of Mathematics\\ University of North Carolina at
Charlotte\\ Charlotte, NC28223} \email{ydiao@uncc.edu}

\author{Claus Ernst}
\address{Department of Mathematics and Computer Science,
Western Kentucky University, Bowling Green, KY 42101}
\email{claus.ernst@wku.edu}

\author{Attila Por}
\address{Department of Mathematics and Computer Science,
Western Kentucky University, Bowling Green, KY 42101}
\email{Attila.por@wku.edu}

\author{Uta Ziegler}
\address{Department of Mathematics and Computer Science,
Western Kentucky University, Bowling Green, KY 42101}
\email{Uta.Ziegler@wku.edu}


\thanks{Y. Diao is currently supported in part by NSF grant DMS-0712958, C. Ernst and U. Ziegler are currently supported in part by NSF grant DMS-0712997.}

%    General info
%\date{January 1, 1994 and, in revised form, June 22, 1994.}


\begin{abstract}
For a knot or link $\K$, let $L(\K)$ be the ropelength of $\K$ and
$Cr(\K)$ be the crossing number of $\K$. In this paper, we show that
there exists a constant $a>0$ such that $L(\K)\le a Cr(\K) \ln^5 (Cr(\K))$ for any $\K$. This result shows that the upper bound of the ropelength  of any knot is almost linear in terms of its minimum crossing number, and is a significant improvement over the best known upper bound established previously, which is of the form $L(\K)\le O(Cr(\K)^{\frac{3}{2}})$.  The approach used to establish this result is in fact more general. In fact, we prove that any 4-regular plane graph of $n$ vertices can be embedded into the cubic lattice with an embedding length at most of the order $O(n\ln^5(n))$, while preserving its topology. Since a knot diagram can be treated as a 4-regular plane graph. More specifically,  Although the main idea in the proof uses a divide-and-conquer technique, the task is highly non-trivial because the topology of the knot (or of the graph) must be preserved by the embedding.
%    Set fboxsep to 0 so that the actual size of the box will match
\end{abstract}

\maketitle

\section{Introduction}\label{sec1}

In the last 3 decades, knot theory has found many important applications in biology \cite{ES,RCV,SECS}. More often than not, in such applications, a knot can no longer be treated as a volumeless simple closed curve in ${\R}^3$ as in classical knot theory. Instead, it has to be treated as a rope like object that has a volume. For example, it has been reported that various knots occur in circular DNA extracted from bacteriophage heads with high concentration and it has been proposed that these (physical) knots can be used as a probe to investigate how DNA is packed (folded) inside a cell \cite{AT,AV1,AV}.
Such applications motivate the study of thick knots, namely knots realized as closed (uniform) ropes of unit thickness.
An essential issue here is to relate the length of a rope (with unit thickness) to those knots that can be tied with this rope.

\medskip
To define the ropelength of a knot, one has to define the thickness of the knot first. There are different ways to define the thickness of a knot, see for example \cite{CKS2,DER3,simonthick}. In this paper, we use the so called {\it disk thickness} introduced in \cite{simonthick} and described as follows. Let $K$ be a $C^2$ knot. A number $r>0$ is said to be {\it nice} if for any distinct points $x$, $y$ on $K$, we have $D(x,r)\cap D(y,r)=\emptyset$, where  $D(x,r)$ and $D(y,r)$ are the discs of radius $r$ centered at $x$ and $y$ which are normal to $K$. The {\it disk thickness} of $K$ is defined to be $t(K)=\sup\{r: r$ is nice$\}$. It is shown in \cite{CKS2} that the disk thickness definition can be extended to all $C^{1,1}$ curves. Therefore, we restrict our discussions to such curves in this paper. However, the results obtained in this paper also hold for other thickness definitions with a suitable change in the constant coefficient.

\begin{definition}
For any given knot type $\K$, a {\it thick realization} $K$ of $\K$ is a knot $K$ of unit thickness which is of knot type
$\K$. The {\it ropelength} $L(\K)$ of $\K$ is the infimum of the length of $K$ taken  over all thick realizations of $\K$.
\end{definition}

The existence of $L(\K)$ is shown in \cite{CKS2}. The main goal of this paper is to establish an upper bound
on $L(\K)$ in terms of $Cr(\K)$, the minimum crossing number of $\K$.

\medskip
It is shown in \cite{Buck,Buck2} that there is a
constant $a>0$ such that for any $\K$, $L(\K)\ge a\cdot
(Cr(\K))^{3/4}$. This lower bound is called the {\em three-fourth power law}. This three-fourth power law is shown to be achievable for some knot families in \cite{CKS1,DE1}. That is, there exists a family of (infinitely many) knots $\{\K_n\}$ and a constant $a_0>0$ such that $Cr(\K_n)\to
\infty$ as $n\to \infty$ and $L(\K_n)\le a_0\cdot (Cr(\K_n))^{3/4}$. On the other hand, it is known that the three-fourth power law does not hold as the upper bound of ropelengths. In fact, it is shown in \cite{DET} that there exists a family of infinitely many prime knots $\{K_n\}$
such that $Cr(\K_n)\to \infty$ as $n\to \infty$ and $L(\K_n)=O(Cr(\K_n))$. That is, the general upper bound of $L(\K)$ in terms of $Cr(\K)$ is at least of the order $O(Cr(\K))$.

In \cite{DEZ2} it is shown that the upper bound of $L(\K)$ is of the order $O(Cr(\K))$ for any Conway algebraic knot. The family of Conway algebraic knots is a very large knot family that includes all 2-bridge knots and Montesinos knots as well as many other knots. The approach used in \cite{DEZ2} is in fact a simpler version of the divide-and-conquer techniques used in this article. The best known general upper bound of $L(\K)$ until now is of order $O((Cr(\K))^{\frac{3}{2}})$, which was obtained in \cite{DEY}. It remains an open question whether $O(Cr(\K))$ is the general ropelength upper bound for any knot $\K$.

\medskip
In this paper, we prove that the general upper bound of $L(\K)$ is almost linear in terms of $Cr(\K)$. More specifically, it is established that there exists a positive constant $a$ such that $L(\K)\le a Cr(\K) \ln^{5}(Cr(\K))$ for any knot $\K$. This is accomplished by showing that a minimum projection of $\K$ can be embedded in the cubic lattice as a planar graph in such a way that the total length of the embedding is of the order at most $O(Cr(\K) \ln^{5}(Cr(\K)))$ and that the original knot can be recovered by some local modifications to this embedding without significantly increasing the total length of this embedding. The construction of the embedding heavily relies on a divide-and-conquer technique that is based on separator theorems for planar graphs \cite{mil}, but many new concepts and results are also needed in this quest. Overall, this is a rather complicated (at least technically) task that requires attention to many technical details.

\medskip
The rest of this paper is organized as follows. In Section \ref{sec2} we introduce the basic concepts of topological plane graphs, cycle cuts and vertex cuts of plane graphs, and some important graph theoretic results concerning these cuts (separator theorems of planar graphs due to Miller \cite{mil}). In Section \ref{sec3}, we introduce the concept of a class of special plane graphs called BRT-graphs. These are the plane graphs we use as basic building blocks when we subdivide knot diagrams. In Section \ref{sec5}, we apply the separator theorems to show the existence of subdivisions of BRT-graphs. It is a typical divide-and-conquer technique that such subdivisions must be ``balanced", that is, each subdivsion produces two smaller BRT-graphs of roughly equal size. We then show that we can apply these concepts recursively subdividing a BRT-graph into smaller and smaller graphs. In Section \ref{sec6}, we introduce the concepts of two special kinds of plane graph embeddings (called ``standard 3D-embedding" and ``grid-like embedding"). These are not lattice embeddings. The grid-like embedding is almost on the lattice and is used as a basic building block to reconstruct the subdivided graphs. The standard 3D-embedding is used as bench mark for verifying the topology preservation of the reconstructed graphs obtained using the grid-like embedding. Section \ref{sec77} is devoted to providing detailed descriptions on how to obtain a grid-like embedding of a plane graph either directly, or indirectly from reconnecting two grid-like embeddings of smaller BRT-graphs obtained in the subdivision process. Then in Section \ref{sec73}, we show that a grid-like embedding obtained from reconstruction using our algorithm preserves the topology of the original graph. Section \ref{sec9} establishes the  upper bound on the length of the embedding generated by our embedding algorithm, from which our main theorem result follows trivially. Finally, we end the paper with some remarks and open questions in Section \ref{sec91}.

\section{Basic Terminology on Topological Plane Graphs and Cycle Cuts of Weighted Plane Graphs}\label{sec2}

Throughout this paper, we use the concept of topologically equivalent graphs. In this case, the vertices are points in $\R^3$ and the edges are space curves that can be assumed to be piecewise smooth. If two edges are incident at a vertex or two vertices, then they intersect each other at these vertices, but they do not intersect each other otherwise. A plane ambient isotopy is defined as a homeomorphism $\Psi:\ \R^2\times [0,1]\rightarrow \R^2$ such that $\Psi(\cdot,t)$ is a homeomorphism from $\R^2$ to $\R^2$ for each fixed $t$ with $\Psi(\cdot,0)=id$.

\medskip
A plane graph $G$ refers to a particular drawing of a planar graph on the plane (with the above mentioned conditions, of course). Two plane graphs $G_1$ and $G_2$ are said to be {\em topologically equivalent} if there exists a plane isotopy $\Psi$ such that $\Psi(G_1,1)=G_2$. It is possible that there are plane graphs $G_2$ that are isomorphic to $G_1$ as graphs, but not topologically equivalent to $G_1$. Since the plane graphs of interest in this paper arise as knot or link projections with the over/under strand information ignored at the crossings, only topologically equivalent plane graphs are considered. At each vertex $v$ of a plane graph $G$, a small circle $C$ centered at $v$ is drawn such that each edge of $G$ incident to $v$ intersects $C$ once (unless the edge is a loop edge in such case the edge will intersect $C$ twice). If $C$ is assigned the counterclockwise orientation, then the cyclic order of the intersection points of the edges incident to $v$ with $C$ following this orientation is called the {\em cyclic edge-order} at $v$. The following lemma assures that two isomorphic plane graphs are topologically equivalent if the cyclic edge-order is preserved by the graph isomorphism at every vertex. This fact can be easily established by induction on the order of the graph. We leave its proof to our reader as an exercise.

\begin{lemma}
\label{equivalence}
Let $G_1$, $G_2$ be two isomorphic plane graphs with $\phi:\ G_1\longrightarrow G_2$ being the isomorphism. If for each vertex
$v$ of $G_1$, the cyclic edge-order of all edges $e_1$, $e_2$, ..., $e_j$ (that are incident to $v$) around $v$ is identical to the cyclic edge-order of $\phi(e_1)$, $\phi(e_2)$, ..., $\phi(e_j)$ around $\phi(v)$, then there exists a plane isotopy $\Psi:\ \R^2\times [0,1]\longrightarrow \R^2$ such that $\Psi(G_1,1)=\phi(G_1)=G_2$. In other words, $G_1$ and $G_2$ are topologically equivalent plane graphs.
\end{lemma}

Frequently, a plane graph needs to be redrawn differently while keeping its topology. These redrawn graphs occur in rectangular boxes and throughout the paper only rectangles and rectangular boxes whose sides are parallel to the coordinate axes are used. It is understood from now on that whenever a rectangle or a rectangular box is mentioned, it is one with such a property. The following simple lemma is also needed later. It can be proven using induction and the proof is again left to the reader.

\begin{lemma}
\label{ncurvesinrec}
Let $R$ be a rectangle and $x_1$, $x_2$, ..., $x_n$, $y_1$, $y_2$, ..., $y_n$ be $2n$ distinct points in $R$, then there exist $n$ disjoint (piecewise smooth) curves $\tau_1$, $\tau_2$, ..., $\tau_n$ such that $\tau_j$ starts at $x_j$ and ends at $y_j$. Furthermore, for any given simply connected region $\Omega$ in the interior of $R$ that does not contain any of the points $x_1$, $x_2$, ..., $x_n$, $y_1$, $y_2$, ..., $y_n$, the curves $\tau_1$, $\tau_2$, ..., $\tau_n$ can be chosen so that they do not intersect $\Omega$. In fact,  many such regions may exist, so long as they do not intersect each other.
\end{lemma}

Using this simple fact, it is possible to redraw any plane graph $G$ in any given rectangle such that the vertices of $G$ are moved to a set of pre-determined points in $R$. This is stated in the following lemma.

\begin{lemma}\label{redrawGinR}
Let $G$ be a plane graph with vertices $v_1$, $v_2$, ..., $v_n$ and let $R$ be a rectangle disjoint from $G$ with $n$ distinct points $y_1$, $y_2$, ..., $y_n$ chosen. Then there exists a plane isotopy $\Psi$ such that $\Psi(G,1)$ is contained in $R$ and $\Psi(v_j,1)=y_j$. Furthermore, for any given simply connected region $\Omega$ in the interior of $R$ that does not contain any of the points $y_j$, $\Psi$ can be chosen so that it keeps $\Omega$ fixed.
\end{lemma}

\begin{proof}
We give a proof for the case that $\Omega= \emptyset$. The case when $\Omega\ne \emptyset$ is left to the reader. A shrinking isotopy $\Psi_1$ is used such that $\Psi_1(G,1)$ is contained in a small rectangle $R_1$ that is small enough to be contained in $R$. Then $R_1$ is moved through a translation to within $R$ such that the vertices $x_j$ of the resulting graph do not overlap with the $y_j$'s. By Lemma \ref{ncurvesinrec}, $x_j$ can be connected to $y_j$ with a curve $\tau_j$ such that $\tau_1$, $\tau_2$, ..., $\tau_n$ do not intersect each other. $\Psi$ can then be obtained by deforming the plane within $R$ by pushing $x_j$ to $y_j$ along $\tau_j$ while keeping the other $\tau_i$ fixed, one at a time.
\end{proof}

\medskip
The following lemma is similar to the above under a different and more restrictive setting. Again this can be proven easily by induction and the proof is left to the reader.

\begin{lemma}\label{fixedcurvelemma}
Let $G$ be a plane graph drawn in a rectangle $R$. Suppose that $R$ contains $n$ disjoint curves $\gamma_1$, $\gamma_2$, ... , $\gamma_n$ which are not closed and are without self intersections. Moreover, these curves do not intersect the vertices of $G$. Let $v_1$, $v_2$, ..., $v_k$ be any $k$ vertices of $G$ ($k\le |G|$) and $x_1$, $x_2$, ..., $x_k$ be any $k$ distinct points in $R$ that are not contained in the curves $\gamma_j$, then there exists a plane isotopy that is identity outside a small neighborhood of $R$ as well as on the $\gamma_j$ and that takes $v_j$ to $x_j$.
\end{lemma}

\medskip
The following is a 3D variation of Lemma \ref{redrawGinR} which is needed later.

\begin{lemma}\label{straighten}
Let $G$ be a plane graph with vertices $v_1$, $v_2$, ..., $v_n$ contained in a rectangle $R\times \{0\}$ in the plane $z=0$. Let  $y_1$, $y_2$, ..., $y_n$ be $n$ distinct points in the rectangle $R\times \{t\}$ in the plane $z=t$ for some number $t>0$. Assume that $y_j$ is connected to $v_j$ by a curve $\nu_j$ that is strictly increasing in the $z$-direction such that the curves $\nu_j$ are disjoint, then there exists a 3D isotopy $\Psi$ such that $\Psi$ is level preserving in the $z$-direction, the identity in the space $z\ge t$ and  outside a small neighborhood of $R\times [0,t]$, and $\Psi(\nu_j,1)$ is a vertical line segment ending in $y_j$ for each $j$.
\end{lemma}

\begin{proof}
This is obvious if there is only one curve $\nu_1$. Assume that this is true for $n=k\ge 1$. Then for $n=k+1$, apply such a level preserving isotopy $\Psi_1$ to the first $k$ curves. $\Psi_1(\nu_{k+1},1)$ is still a strictly increasing simple curve from $\Psi_1(v_{k+1},1)$ to $y_{k+1}$. Modify $\Psi_1(\nu_{k+1},1)$ in a level preserving fashion so that its projection to the plane $z=0$ is a simple curve without self intersection. Now a plane isotopy can be constructed by a push back along this curve from $\Psi_1(v_{k+1},1)$ to the projection $y^\p_{k+1}$ of $y_{k+1}$ to the plane $z=0$. This does not affect the projections $y^\p_j$ of the points $y_i$ into the plane $z=0$ for $j\le k$. This plane isotopy can be used to define the level preserving isotopy that works for $n=k+1$ curves. Although there are still some technical details in the argument, it is intuitively obvious at this point and the details are left to the interested reader.
\end{proof}

\begin{remark}\label{re22}
If we relax the condition that the simple curves $\nu_j$'s are strictly increasing in the $z$ direction to that they are non-decreasing in the $z$ direction, then one can show that the result of Lemma \ref{straighten} holds without the requirement that the isotopy is level preserving. It is important to note that in this case $\Psi(\cdot,0)$ is still the identity and $\Psi(G,1)$ is topologically equivalent to $G$.
\end{remark}

\medskip
A 3-dimensional ambient isotopy can be similarly defined on $\R^3\times [0,1]$. However we have to be much more careful about the cyclic edge order at the vertices. For a plane graph the cyclic order of the edges at a vertex is well defined however for a graph in 3 dimensional space it is not. We will fix this problem by requiring that for a graph in 3 dimensional space each vertex $v$ is contained in a small 2 dimensional disk-neighborhood $D_v$. All the edges that terminate at $v$ intersect $D_v$ in a short arc that terminates at $v$. The intersection points of these arcs on $D_v$ now define the cyclic edge order at $v$. We will require that all 3-dimensional isotopies preserve this structure, that is, a 3-dimensional isotopy can move the small disk $D_v$ around in 3-space, it can even deform it, however throughout the isotopy $v$ remains on $D_v$ and all edges terminating at $v$ keep their short arcs on $D_v$. In this way the cyclic edge order at $v$ remains invariant under the 3-dimensional isotopy. As we will see in Section \ref{sec6}, we use two types of such neighborhoods called blue triangles and red square.  Since all 3D-isotopies used in this paper are ambient isotopies that preserve the neighborhood structure of a vertex, we will call them {\em 3D VNP-isotopies} (where VNP stands for ``Vertex Neighborhood Preserving").

\medskip
\begin{definition}
Let $G$ be a plane graph. $G$ is called a {\em weighted graph} if each vertex, edge, and face of $G$
is assigned a weight (i.e., a non-negative number) and the sum of these weights is 1.
\end{definition}

\begin{definition}
Let $G$ be a weighted plane graph. A {\em cycle cut} of $G$ is a cycle $\gamma$ in $G$ such
that deleting all vertices in $\gamma$ (and the edges connected to them) divides $G$ into two subgraphs $G_1$, $G_2$, on opposite sides of $\gamma$, i.e. $G\setminus \gamma = G_1\cup G_2$.
Moreover, if the weights in each $G_i$ sum to no more than $\alpha$ for some real number
$\alpha$, $0<\alpha <1$, then the cycle cut is called an $\alpha$-cycle cut. The {\it size} of the cycle cut $\gamma$ is the number of vertices in $\gamma$.
\end{definition}

The following theorems are proved in \cite{mil} and play vital roles in the proof of the main theorem of this paper.

\begin{theorem}
\label{cyclecut}
\cite{mil}
Let $G$ be a 2-connected and weighted plane graph such that no face of $G$ has weight
more than $2/3$, then there exists a $\frac{2}{3}$-cycle cut of size at most
$2\sqrt{2\lfloor d/2 \rfloor n}$, where $d$ is the maximal face size of $G$.
\end{theorem}

If one drops the assumption that $G$ is 2-connected then Theorem \ref{cyclecut} may not hold since $G$ may be a tree. In this case there is the following theorem.

\begin{theorem}
\label{generalcut}
\cite{mil}
Let $G$ be a connected and weighted plane graph such that all faces have been assigned weight zero, then either there exists a $\frac{2}{3}$-cycle cut of size
$2\sqrt{2\lfloor d/2 \rfloor n}$, or there exists a cut vertex $v$ of $G$ such
that each connected component in $G\setminus v$  has a total weight less or equal to $2/3$.
\end{theorem}

\medskip
We refer to the cycle $\gamma$ that generates a cycle cut as defined by Theorems \ref{cyclecut} and \ref{generalcut} as a {\em cut-cycle}.

\section{Component-Wise Triangulated Plane Graphs: Definitions and Subdivisions}\label{sec3}

A main tool used to achieve the ropelength upper bound obtained in this paper is the ``divide-and-conquer" approach familiar to researchers in graph theory. That is, a knot projection (treated as a 4-regular plane graph) is divided repeatedly using the theorems given in the last section. Several non trivial issues arise in this process. First, the graphs obtained after repeated subdivisions may become highly disconnected. Second, one needs to keep track of the topology of the subdivided piece of the graph so that the pieces resulting from the subdivisions can later be assembled correctly and the original graph can be recovered. This section describes in detail how we handle these problems. The aim is to impose a special structure on the graphs that is preserved after each cut and that this structure allows us to reconstruct the original graph (with the correct topology) from the pieces generated in the divide-and-conquer process.

\begin{definition}
Let $G$ be a connected plane graph without loops. We say that $G$ admits a {\em proper BR-partition} if the vertices of $G$ can be partitioned into two sets $V^B$ and $V^R$ (called {\it blue} and {\it red} vertices, respectively) such that there is no edge between any two blue vertices.
\end{definition}

Let $G$ be a plane graph that admits a proper BR-partition with $V^B$ and $V^R$ being the set of blue and red vertices respectively. $V^R$ induces a subgraph $G(V^R)$ of $G$ that is itself a plane graph. A connected component of $G(V^R)$ is called a {\it red component}. For a red component $M$,
let $V^R_M$ be its vertex set and let $V^B_M$ be the set of blue vertices that are adjacent to some vertices in $V^R_M$. Let $V^*_M=V^R_M \cup V^B_M$ and let $G(V^*_M)$ be the (plane) subgraph of $G$ induced by $V^*_M$. $G(V^*_M)$ is called a {\it BR-component} of $G$. Notice that under this definition, different BR-components may share common blue vertices, but each edge $e$ of $G$ belongs to exactly one BR-component. A graph is said to be {\em triangulated} if each face of the graph is either a triangle or a digon.

\begin{definition}
\label{CWdef}
Let $G$ be a connected plane graph that admits a proper BR-partition. $G$ is called a {\it BRT-graph} if each BR-component $G(V^*_M)$ of $G$ is
triangulated. A triangulated BR-component $G(V^*_M)$ is called a {\it BRT-component}.
\end{definition}

\begin{lemma}\label{L3.3}
Let $G$ be a BRT-graph, then (i) the boundary of any face of $G$ contains at most one blue vertex and (ii) any cut vertex of $G$ is a blue vertex.
\end{lemma}

Notice that condition (ii) above is equivalent to the following statement: for any red component $M$,
the BRT-component $G(V^*_M)$ is 2-connected.

\medskip
\begin{proof}
Suppose that the boundary $\partial F$ of a face $F$ in $G$ contains two different blue vertices $v$ and $w$. Let $P$ be a
path in $\partial F$ that connects $v$ and $w$ which contains at least three vertices, since $P$ cannot be a single blue-blue edge. We assume that $v$ and
$w$ and $P$ were chosen so that there are only red vertices on $P$ besides $v$ and $w$. $v$ and $w$ must belong to $G(V^*_M)$ where $M$ is the red component containing the red vertices on $P$. It is easy to see that even if $P$ contains only a single red vertex, $\partial F$ must contain at least 4 edges since $v$ and $w$ cannot be connected by a single edge. This implies that $G(V^*_M)$ is not triangulated, which contradicts the given condition that $G$ is a BRT-graph. This proves (i).

\medskip
Let $v$ be a cut vertex of $G$. Then there exists a face $F$ of $G$ such that $v$ appears on $\partial F$ at least twice, i.e., $\partial F$ can be described by a walk $w=P_1P_2$ where $P_1$ and $P_2$ are closed walks in $G$ starting and ending at $v$. $P_1$ cannot contain only one edge since that would force a loop edge. Thus $P_1$ contains at least two edges. Similarly, $P_2$ also contains at least two edges. Thus there exist vertices $w_1$ and $w_2$ on $P_1$ and $P_2$, respectively, that are adjacent to $v$. Note that $w_1w_2$ cannot be an edge in $G$ since $v$ is a cut-vertex.
Let us assume that $v$ is a red vertex, and let $M$ be the red component that contains $v$. Then $w_1$ and $w_2$ must also be contained in $G(V^*_M)$ and must be on the same face of $G(V^*_M)$. However this contradicts the assumption
that $G(V^*_M)$ is triangulated. Thus $v$ must be blue. This proves (ii).
\end{proof}

For a BRT-graph $G$, let us define a graph $T_G$ that is associated with $G$ in the following way. Each blue vertex of $G$ corresponds to a vertex in $T_G$ (which is still called a blue vertex) and each red component of $G(V^R)$ corresponds to a vertex in $T_G$ (which is called a red vertex). Vertices of the same color in $T_G$ are never adjacent and a blue vertex $x$ and a red vertex $y$ in $T_G$ are connected by a single edge if and only if the blue vertex $v$ in $G$ corresponding to $x$ is contained in $G(V^*_M)$ where $M$ is the red component corresponding to $y$. The following lemma asserts that the graph $T_G$ so constructed is a tree.

\begin{lemma}\label{CycleinBRT}
The graph $T_G$ defined above is a tree. Equivalently, any cycle in $G$ is contained in a single BRT-component.
\end{lemma}

\begin{proof}
By construction the graph $T_G$ is a simple bipartite graph since two vertices of different colors can be connected by at most one edge and vertices of the same color are never connected. Moreover there can be no digon in $T_G$ by construction.
It follows that a cycle $C$ in $T_G$ (if it exists) must contain at least 4 vertices and the colors of the vertices on the cycle must alternate if one travels along $C$. Assume that there is a cycle $C$ in the graph $T_G$ which contains a path $y_1, x_1, y_2, x_2, y_3$ where the $x_i$ are red and the $y_i$ are blue. Let $M$ be the red component of $G(V^*_M)$ which corresponds to $x_1$. Let $y_1$ and $y_2$ be the two blue vertices in $C$ that are adjacent to $x_1$. Let $v_1$ and $v_2$ be the two blue vertices in $G$ that correspond to $y_1$ and $y_2$, respectively. The fact that $y_1$ and $y_2$ are connected to $x_1$ in $T_G$ implies that $v_1$ and $v_2$ are adjacent to some vertices in $M$.

The red component $M_1$ of $G(V^R)$ that corresponds to $x_1$ divides the plane into one outer face $F$ and several (possibly zero) inner faces. In the interior of each face of $M_1$ there are blue vertices of $G(V^*_M)$ (if any) or vertices of the other BRT-components besides $M_1$ (if any).
If $v_1$ and $v_2$ are both contained in the interior of $F$, then there exists a path from $v_1$ to $v_2$ of the form $v_1r_1...r_kv_2$ where $k\ge 1$ and the path $r_1...r_k$ is contained in $\partial F$. Since the face in $G(V^*_{M_1})$ containing this path cannot be triangulated without using blue-blue edges, this is not possible. Thus $v_1$ and $v_2$ cannot be both contained in the interior of $F$. Similarly, on can show that $v_1$ and $v_2$ cannot be both contained in the interior of the same face for any face of $M_1$. Therefore, one of them, say $v_1$ is contained in the interior of an inner face $F_1$ of $M_1$ and the other $v_2$ is either contained in the interior of a face $F_2$ which is either the outer face $F$ or a different inner face.

\medskip
Since the cycle $C$ has at least length four, there exists a second red vertex $x_2$ in $C$ neighboring $y_2$. Let $M_2$ be the red component of $G(V^R)$ corresponding to $x_2$. This means that $v_2$ is connected by an edge to a red vertex in $M_2$. Since $v_2$ is contained in the interior of $F_2$ and $G$ is a plane graph, this implies that $M_2$ is entirely contained in the interior of $F_2$. Similarly, continuing to move along the cycle $C$ in the direction established by moving from $y_2$ to $x_2$, leads to the conclusion that all blue vertices and red components (of $G$) corresponding to the vertices on $C$ are contained in the interior of $F_2$. This contradicts the fact that $v_1$ is not contained in the interior of $F_2$. Thus the cycle $C$ cannot exist.

\medskip
It is easy to see that any cycle passing through more than one
BRT-component gives rise to a cycle in $T_G$. Thus any cycle in $G$ is contained in a single BRT-component of $G$.
\end{proof}

Let $G$ be a BRT-graph. Let us consider two different ways of dividing $G$ into subgraphs. First, consider the case that $G$ has a cut-vertex $v$ (recall that $v$ must be a blue vertex by Lemma \ref{L3.3}). Let $G_1$, $G_2$, ..., $G_k$ be the
connected components of $G\setminus \{v\}$ with $V_i$ being the set of vertices of $G_i$.

Pick an arbitrary proper subset $I_1$ of $I=\{1,2,...,k\}$ and let $I_2 = I \setminus I_1$.
Let $U_1 = \cup_{i\in I_1} V_i \cup \{v\}$ and $U_2 = \cup_{j\in I_2} V_j \cup \{v\}$. Let $J_1$ be the induced subgraph of $G$ consisting of the BRT-components whose vertices are contained in $U_1$ and $J_2$ be the induced subgraph of $G$ consisting of the BRT-components whose vertices are contained in $U_2$. Note that $J_1$ and $J_2$ inherit the plane graph structure naturally from $G$. In particular,
the cyclic ordering of the edges around the vertex $v$ in each subgraph $J_1$, $J_2$ is naturally inherited from the cyclic edges ordering around $v$ in $G$ This describes the first kind of subdivision of $G$ which is formally defined below.

\begin{definition}
\label{subdivisionvertex}
Let $G$ be a BRT-graph with a cut-vertex $v$, then dividing $G$ into two subgraphs $J_1$ and $J_2$ as described above is called {\em subdividing $G$ by a vertex-cut}.
\end{definition}

Before describing the second subdivision which depends on a cycle of $G$, we introduce the concept of a {\em normal cycle}. Let $u$, $v$, and, $w$ be three vertices of $G$ such that $uvwu$ is a triangle (namely a cycle with three edges) in $G$.
T he triangle $uvwu$ is empty if it is the boundary of a face of $G$. A cycle $\gamma$ in $G$ is said to be {\it normal} if it contains at least two vertices and no three consecutive red vertices on $\gamma$ form an empty triangle in $G$.

\medskip
\begin{lemma}
\label{cutcycleexistance}
Let $G$ be a BRT-graph. If $G$ has a cut-cycle, then it has a normal cut-cycle.
\end{lemma}

\begin{proof}
Assume that $P$ is a cut-cycle in $G$ with length $\ell$. Then $\ell$ is at least $2$, since $G$ does not contain any loops. If $P$ contains three consecutive red vertices $u$, $v$, $w$ which form an empty triangle, then remove $v$ and replace $uvw$ by $uw$. Since this operation can never reduce the number of vertices below 2, it always leads to a normal cut-cycle with length at least 2.
\end{proof}

\begin{figure}[htbp]
\begin{center}
\includegraphics[scale=0.4]{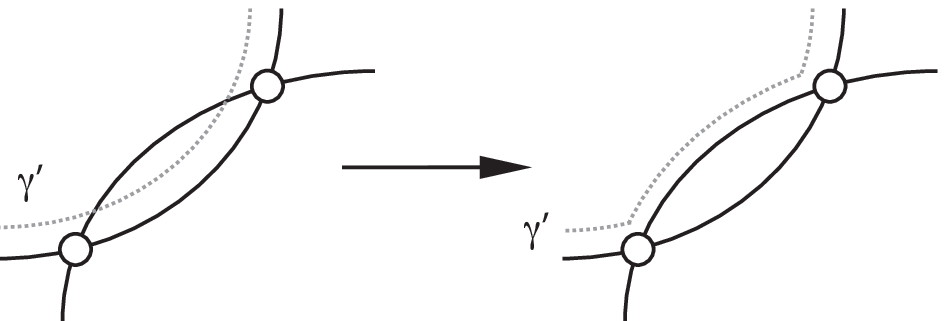}
\hskip 0.3in
\includegraphics[scale=0.5]{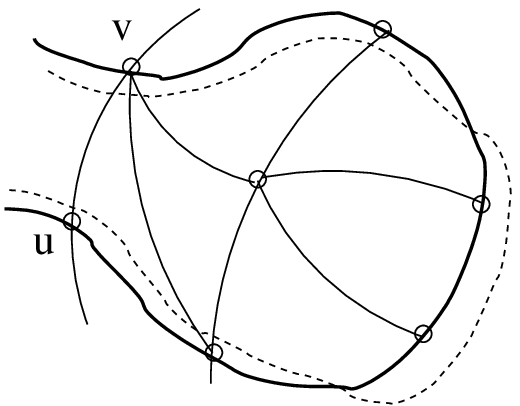}
\hskip 0.3in
\includegraphics[scale=0.5]{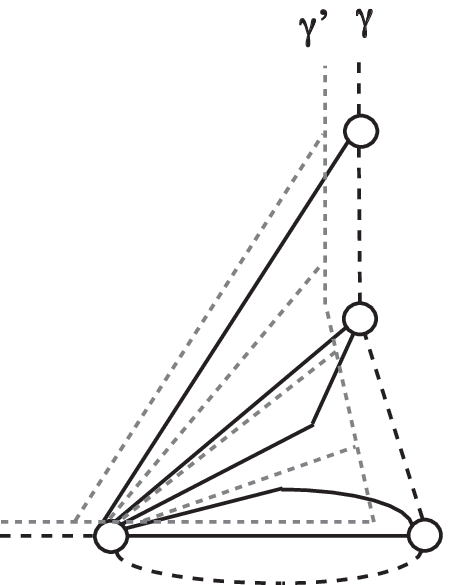}

\caption{$\gamma^\p$ intersects an edge twice. Left: Removing empty digons by re-routing $\gamma^\p$. Middle: An edge intersected by $\gamma^\p$ (dashed) twice in a non-trivial way. Right: Pushing $\gamma^\p$ across several digons. }\label{digonavoid}
\end{center}
\end{figure}

\medskip
Now let us describe how to use a normal cut cycle $\gamma$ to divide $G$. First, $\gamma$ is pushed off the red vertices by a small distance, resulting in a simple closed curve $\gamma^\p$. The push off $\gamma^\p$ is required to have the following properties: (1) All intersections of $\gamma^\p$ with the edges of $G$ must happen transversely (except possibly at their endpoints if they are connected to a blue vertex);
(2) If $v$ is a blue vertex on $\gamma$ then $v$ stays on $\gamma^\p$; (3) If $e$ is an edge on $\gamma$ that is connected to a blue vertex $v$ on $\gamma$, then $\gamma^\p$ does not intersect $e$; (4) If $e$ is an edge on $\gamma$ that is connected to two red vertices, then $\gamma^\p$ may intersect $e$ at most once; (5) For an edge $e$ that is not on $\gamma$ but is connected to exactly one red vertex on $\gamma$, $\gamma^\p$ can intersect $e$ at most once; (6) For an edge $e$ that is not on $\gamma$ but is connected to two red vertices on $\gamma$, $\gamma^\p$ can intersect $e$ at most twice. (The right of Figure \ref{digonavoid} shows a non-trivial example of this.) It is easy to see that such a push off is always possible.
In addition, if $\gamma^\p$ intersects $e$ twice then no empty digons (an empty digon is a digon whose interior or exterior does not intersect $G$) should be created. The empty digons can always be avoided by routing $\gamma^\p$ around the digon as shown in Figure \ref{digonavoid}. Note that such digons be nested and $\gamma^\p$ may have to be pushed across several digons, see Figure \ref{digonavoid} on the right.

A simple closed curve $\gamma^\p$ obtained from a cycle $\gamma$ with these properties is called a {\em push-off} of $\gamma$.

\medskip
Let $\gamma$ be a normal cycle of $G$ and let $\gamma^\p$ be a push-off of $\gamma$. Let us (temporarily) insert into $G$ a white colored vertex at each intersection point of $\gamma^\p$ with the edges of $G$ and a new edge for every arc on $\gamma^\p$ connecting two such white vertices or connecting a white vertex and a blue vertex. These white vertices are considered the vertices of $\gamma^\p$ and are denoted by $V(\gamma^\p)$. The arcs of $\gamma^\p \setminus V(\gamma^\p)$ are considered as the edges of $\gamma^\p$.

\begin{figure}[htbp]
\begin{center}
\includegraphics[scale=0.6]{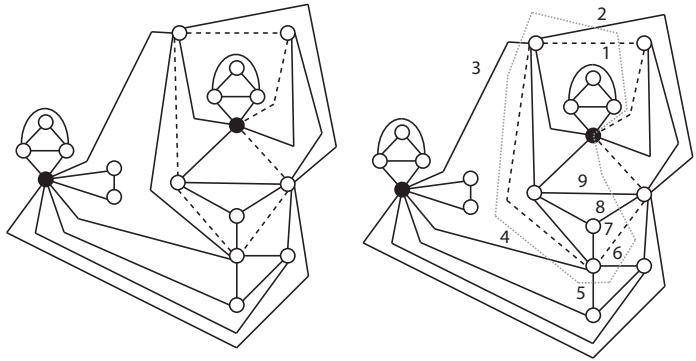}
\caption{ Subdividing a plane graph by a edge-cut. Left: a cut cycle $\gamma$ shown by the dashed edges. Right: the modified normal $\gamma$ (dashed) and its push-off $\gamma^\p$ (dashed and grey). Blue vertices are marked by dark circles and red vertices are marked by while circles.
The edges cut by $\gamma^\p$ are labeled in cyclic order.}
\label{Gcyclecut}
\end{center}
\end{figure}

\begin{figure}[htbp]
\begin{center}
\includegraphics[scale=0.6]{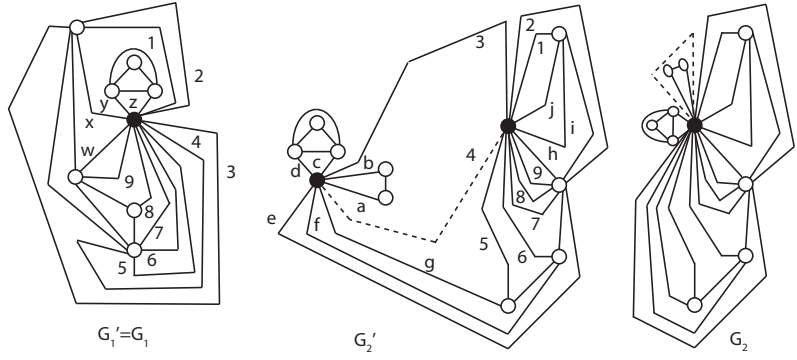}
\caption{ The new graphs obtained by subdividing the graph in Figure \ref{Gcyclecut} by $\gamma$. Left: $G_1$ obtained from the subgraph inside $\gamma^\p$. Middle: the graph $G^\p_2$ obtained from the subgraph outside $\gamma^\p$. Right: $G_2$, obtained by contracting edge 4 in $G^\p_2$. The loop edge created by this contraction (dashed) is deleted from $G_2$.}\label{Gsubdivided}
\end{center}
\end{figure}

\begin{figure}[htbp]
\begin{center}
\includegraphics[scale=0.6]{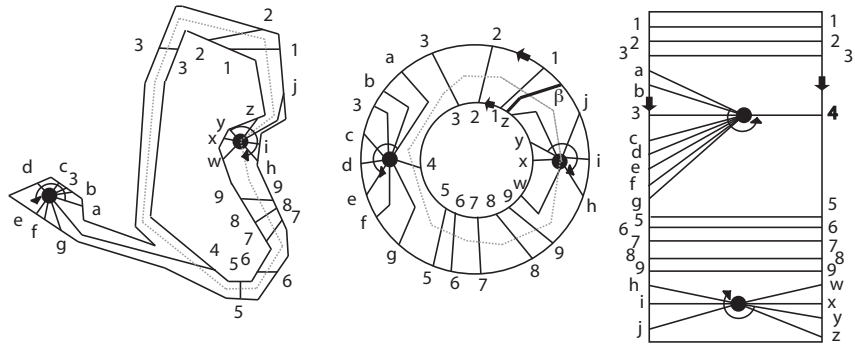}
\caption{Left: the neighborhood $N$ of the edge-cut subdivision of Figures \ref{Gcyclecut} and \ref{Gsubdivided} as it arises in $G$;
Center: $N$ deformed into an annulus together with the path $\gamma^\p$, the orientations on $\partial N$ and the path $\beta$ used to establish the linear order; Right: $N$ (cut open along $\beta$) deformed into a rectangle.}
\label{unicycleneighborhood}
\end{center}
\end{figure}

This graph obtained by adding the edges and vertices of $\gamma^\p$ is called $G^w$. Let $G^w\setminus \gamma^\p$ denote the plane graph obtained from $G^w$ by deleting all vertices $V(\gamma^\p)$ and the edges connected to these vertices. $G^w\setminus \gamma^\p$ is separated into two disjoint subgraphs $G_I$ and $G_O$ with $G_I$ inside of  $\gamma^\p$ and $G_O$ outside of $\gamma^\p$. Let $G_I^*$ be the subgraph of $G^w$ induced by the vertices $V(G_I)\cup V(\gamma^\p)$. Similarly, let $G_O^*$ be the subgraph of $G^w$ induced by the vertices $V(G_O)\cup V(\gamma^\p)$. Finally, in $G_I^*$ we contract $\gamma^\p$ to a single vertex and mark it as a new blue vertex $v_1$. The resulting plane graph is denoted by $G^\p_1$. Similarly, in $G_O^*$ we contract $\gamma^\p$ to a single vertex and mark it as a new blue vertex $v_2$ and call the resulting plane graph $G^\p_2$. See Figures \ref{Gcyclecut}, \ref{Gsubdivided}, and \ref{unicycleneighborhood} for an illustration of this process.

Notice that the cyclic orders of the edges connected to the newly created blue vertices $v_1$ in $G^\p_1$ and $v_2$ in $G^\p_2$ are inherited from the cyclic order of the intersection points of $G$ with $\gamma^\p$.

Thus $G_I^*$ and $G_O^*$ can be recovered from $G^\p_1$ and $G^\p_2$ by expanding $v_1$ or $v_2$ back to $\gamma^\p$.
$G$ can be recreated by gluing $G_I^*$ and $G_O^*$ along $\gamma^\p$ with the original cyclic order of the edges along $\gamma^\p$ preserved.

\medskip
Notice that $G^\p_1$ and $G^\p_2$ may not admit a proper BR-partition for two reasons. First, $G^\p_1$ and $G^\p_2$ may contain loop edges. Second, the newly created blue vertices $v_1$ and $v_2$ (from the contraction of $\gamma^\p$) may be adjacent to some other blue vertices in $G^\p_1$ or $G^\p_2$ which already exist in $G$. If this happens in $G_i$, one of the edges connecting these blue vertices to $v_i$ is simply contracted.

More precisely, assume that $w$ is a blue vertex in $G^\p_i$ that is connected by several edges $e_1,\ldots,e_k$ to $v_i$. One of these edges is picked, say $e_1$ and contracted. This combines $w$ and $v_i$ into one blue vertex - still denoted $v_i$, and generates $k-1$ loop edges at $v_i$, see Figure \ref{Gsubdivided}. After all blue-blue edges (that are not loop edges) have been eliminated in this way, the resulting (plane) graphs are denoted by $G_{1l}$ and $G_{2l}$, the $l$ indicating that the graphs may contain loop edges. A loop edge, after it is created, is never cut again, nor does it influence any further subdivisions. Thus there is no reason to keep these loop edges in the plane graphs for the future subdivision process. After deleting
the loop edges from $G_{1l}$ or $G_{2l}$, the resulting graphs are denoted by $G_1$ and $G_2$. Information about $\gamma^\p$, the contracted edges, and the deleted loops which are not included in $G_1$ or $G_2$ are kept as described later (see Definition \ref{glueinginstructions}). Only $G_1$ and $G_2$ are used in the subsequent subdivisions.

\begin{definition}
\label{subdivisioncycle}
Let $G$ be a BRT-graph with a normal cycle $\gamma$, then dividing $G$ into two plane graphs $G_1$ and $G_2$ using a push-off of $\gamma$ as described above is called {\em subdividing $G$ by an edge-cut}.
\end{definition}

After a BRT-graph $G$ is subdivided into two new plane graphs as
described in the above two definitions, are the newly obtained graphs also BRT-graphs? This is not obvious in the case of an edge-cut subdivision since some original BRT-components may have been modified by the subdivision process. The answer to this is affirmative and established in Lemma \ref{splittinglemma} below.

\begin{lemma}
\label{splittinglemma}
Let $G$ be a BRT-graph. If $G_1$ and $G_2$ are the two graphs obtained after a vertex-cut subdivision or an edge-cut subdivision is applied to $G$, then $G_1$ and $G_2$ are also BRT-graphs.
\end{lemma}

\begin{proof}
In the case that $G_1$ and $G_2$ are obtained by a vertex-cut subdivision the lemma is obvious.
Thus we concentrate on the case of an edge-cut subdivision. By construction the new graphs $G_1$ and $G_2$ do not have blue-blue edges. Furthermore, if $G$ is connected, then $G_1$ and $G_2$ are also connected. Let $\gamma$ be the cycle used in this edge-cut subdivision. By Lemma \ref{CycleinBRT}, $\gamma$ is contained in a single BRT-component $G(V^*_M)$ for some red component $M$ in $G$. Obviously, all other BRT-components of $G$ remain unchanged in the subdivision process. These BRT-components remain as triangulated BRT-components in either $G_1$ or $G_2$.

\medskip
Let $N$ be a red component of $G_1$ that contains some vertices of $M$. To show that $G(V^*_N)$ is triangulated, it is first shown that all faces in $G^\p_1$ or $G^\p_2$ created in the process of contracting $\gamma^\p$ are either triangles or digons. Note that $G^\p_1$ or $G^\p_2$ may contain loop edges and those are addressed later.  To show that after contracting $\gamma^\p$ the resulting graph is still triangulated, the contraction of the edges on $\gamma^\p$ is considered one edge at a time.

\medskip
Let $F$ be a face of $G(V^*_M)$ and let $\partial F$ be its boundary. If the interior of $F$ does not intersect $\gamma^\p$, then $F$ is not affected by the contraction process of $\gamma^\p$ and it remains a triangle or a digon after the subdivision.
If the interior of $F$ and $\gamma^\p$ intersect each other, then $\partial F$ must contain at least one vertex on $\gamma$. Let $e$ be an edge of $\gamma^\p$ that intersects $F$ and assume that $e$ splits $F$ into two faces $F_1$ and $F_2$ in $G^*_I$ and $G^*_O$, respectively. $e$ can intersect $\partial F$ in two white vertices (i.e. in two different edges) or $e$ can intersect $\partial F$ in one white  and one blue vertex (i.e. in one edge and the blue vertex of $F$). It is easy to see that regardless of whether $F$ is a digon or a triangle and regardless of the particular location of these intersections the contraction process changes  $F_1$ and $F_2$ into loops, digons, or triangles, see Figure \ref{BRTpreservation}.

\begin{figure}[htbp]
\begin{center}
\includegraphics[scale=0.6]{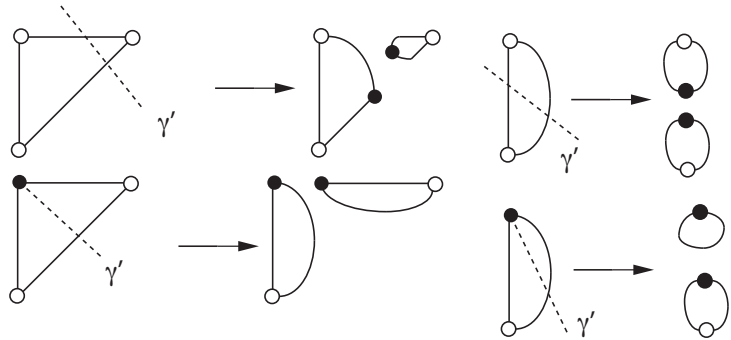}
\caption{The possible cases of how an edge $e$ on $\gamma^\p$ can intersect the interior of a triangle or a digon. The three cases for a triangle are on the left, the two cases for a digon are on the right.  The resulting faces after $\gamma^\p$ (marked by the dashed line) is contracted to the new blue vertex are on the right side of the arrows.}
\label{BRTpreservation}
\end{center}
\end{figure}

\medskip
At this point we have established that after the contraction of $\gamma^\p$ all faces that are changed by the subdivision remain triangulated in $G^\p_1$ and $G^\p_2$. To obtain $G_{1l}$ and $G_{2l}$, it is necessary to contract some of the blue-blue edges that may have been created when the new blue vertices were introduced. However, the contraction of an edge in a graph does not increase the size of any face, and thus the triangulation property of all affected faces is preserved. The final step to obtain $G_1$ and $G_2$ from $G_{1l}$ and $G_{2l}$ is to delete all loop edges that may have been created in the contraction process of $\gamma^\p$ and some of the blue-blue edges. Let $e$ be a loop edge created by this contraction connected to a blue vertex $v$. Let $F_1$ and $F_2$ be the two faces on the different sides of the loop edge $e$. If both $\partial F_1$ and $\partial F_2$ contain vertices different from $v$ then these vertices are red and belong to two different BRT-components (after the cut). Therefore the face $\partial F_1\cup\partial F_2$ created by deleting $e$ does not have to be triangulated and $e$ can be deleted. If one, say $\partial F_1$ contains no other vertices (that is $\partial F_1$ is consists of one or two edges from $v$ to $v$) then the deletion of $e$ causes the number of edges in $\partial F_2$ to remain the same or to decrease by one and thus cannot violate the triangulation property of $F_2$. Therefore all BRT-components in $G_1$ and $G_2$ remain triangulated.
\end{proof}

\medskip
After a BRT-graph $G$ is subdivided into BRT-graphs $G_1$ and $G_2$,
it may no longer be possible to reconstruct $G$ from $G_1$ and $G_2$ unless
information about the subdivision step is kept. For each blue vertex $v$ in $G$ which is not involved in the subdivision process, the order of its edges in $G$ is the order of the edges of the corresponding blue vertex in $G_1$ or $G_2$ and no additional information must be kept. The information about the blue vertices involved in the subdivision process is captured in a small neighborhood $N$. $N$ includes information about the edges involved in the subdivision and their relative order. In addition, for an edge-cut, the neighborhood contains information about the edges contracted (if any) and the loop edges (if any) temporarily created during the subdivision. All loop edges are deleted in the final step of the process that changes $G_{il}$ to $G_i$.

\medskip
Notice that all the edges in $G$ are labeled at the beginning and these labels do not change in the subdivision process. Thus the edges in $G_1$, $G_2$, and $N$ share the same label if and only if they are part of the same edge in $G$, therefore the labeled graphs $G_1$ and $G_2$ together with $N$ contain all the information needed to reconstruct $G$ (since Lemma \ref{equivalence} can then be applied to the reconstructed graph).

\medskip
The detailed information stored in $N$ is different for a circular
edge-cut and for a vertex cut.

\medskip
(i) For a circular edge-cut, $N$ is an annulus which contains a
small neighborhood of $\gamma^\p\cup E_\gamma$, where $E_\gamma$ is the
set of edges of $G^w$ that are incident to both a white and a blue vertex and that are contracted in the steps of the subdivision process which changes $G^\p_{i}$ to $G^\p_{il}$, see Figure \ref{unicycleneighborhood} for an example.

\medskip
The outside and inside boundaries of $N$ are used to keep track of the
cyclic order of the edges of $G_1$ and $G_2$ (as well as the deleted loop edges) around their corresponding blue vertices created by the edge-cut. For reasons discussed in Section \ref{sec6} and in Subsection \ref{sec7} we impose a {\it linear order} on the edges around a blue $v^b_i$ at the time when it is created in the subdivision process. This is accomplished by identifying a path $\beta$ in $N$ that connects the inside and outside boundary of $N$ without intersecting any of the edges of $G$. Cutting $N$ along $\beta$ results in the linear order of the edges for $v^b_1$ and $v^b_2$ inherited from the counterclockwise orientation on each of the boundary components of $N$, see Figure \ref{unicycleneighborhood}. From now on, it is understood that the linear order at a blue vertex is so defined if the vertex is created by a circular edge-cut subdivision of a BRT-graph.

\medskip
A loop may be created and becomes part of $G_{il}$ in two situations. First, if an edge $e$ not on $\gamma$ is cut twice by $\gamma^\p$, (see the right diagram in Figure \ref{digonavoid} for an example)  a loop is created in the $G_{il}$ which does not contain the vertices incident to $e$. In this case the same edge label appears twice on each boundary component of $N$.

Second, if $k$ red-blue edges $e_1,\ldots,e_k$ with $k > 1$ which are incident to a blue vertex $u$ are cut by $\gamma^\p$ (see Figures \ref{Gcyclecut} and \ref{Gsubdivided} for an example), then $u$ is contracted into the $v^b_i$ associated with the $G^\p_i$ which contains $u$ and $k-1$ loops are created.
The loops are in the $G_{il}$ associated with that $v^b_i$. Here the same edge label appears three times in $\partial N$ for each of these loops: twice on the boundary of the component which contains the loop edge before its deletion and once on the other.  Loops are not included in $G_i$. However, to enable a correct reconstruction of $G$ the information about how the loop connections interleave with other edges cut by $\gamma^\p$ is kept by the edge labels on  $\partial N$.

\medskip
(ii) For a vertex-cut subdivision, $N$ is a disk which contains a small neighborhood of $v$, the cut-vertex used for the subdivision. Once a single point $\beta$ that does not belong to any edge on the boundary of $N$ is chosen, the linear order of the edges around the new blue vertices $v_1^b$ and $v_2^b$ is inherited from the cyclic order of $v$.

\begin{figure}[htbp]
\begin{center}
\includegraphics[scale=0.65]{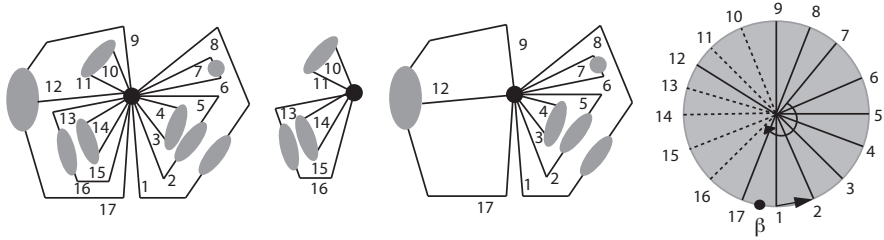}
\caption{Left: The graph $G$ with a vertex-cut at a blue vertex $v$. The gray shaded ellipses are the different BRT-components that contain the vertex $v$. These BRT-components can be of any size; Middle: The graphs $G_1$ and $G_2$ after the vertex-cut subdivision; Right: The neighborhood $N$ of the blue vertex $v$ containing the gluing instruction, the edges of $G_1$ are dashed. The arrow on $\partial N$ shows the orientation and the point $\beta$ is used to define the linear order of the edges.}
\label{vertexcut}
\end{center}
\end{figure}

\medskip
From now on, it is understood that the linear order is so defined at a blue vertex created by a vertex-cut subdivision of a BRT-graph, see Figure \ref{vertexcut}.
Note that we can indicate the linear order of the edges at a blue vertex $v$ in $N$ by a small circular arrow $o_v$ around a blue vertex. The edge on which the tail of $o_v$ is placed indicates the first edge of the linear order and the arrow head points into the direction of that linear order.  Even though we could think of a BRT-graph $G$ as a graph where every blue vertex has an arrow indicating a linear order we are only interested in assigning a linear order (with an arrow) resulting from subdivsions. When we apply an embedding algorithm to a BRT-graphs later in this paper, every blue vertex originated from a subdivision process and has an assigned linear order. The linear order of the blue vertices in a neighborhood $N$ is included in $N$. For an example see the two small arrows around two blue vertices in Figure \ref{unicycleneighborhood}. Note that we did not include the small arrows at the blue vertices in Figures \ref{Gcyclecut} and \ref{Gsubdivided}, since the linear orders were not relevant to our discussions at that time.

\begin{definition}
\label{glueinginstructions}
The small neighborhood $N$ addressed above, together with the labels of the edges of $G_1$,  $G_2$, is called {\em the gluing instruction} of the corresponding subdivision of $G$. The orientation of a component of $\partial N$ is called the {\em orientation vector} of the boundary component.
\end{definition}

The information provided by $G_1$, $G_2$, and the gluing instruction suffices to reconstruct the graph $G$ uniquely (up to a plane isotopy) by observing that the position of the loop edges (up to a plane isotopy) in $G^\p_{il}$ can be derived from the given information. For each such loop edge the gluing instruction determines a unique vertex and a unique face that must contain the loop edge and therefore the position of the loop edge is unique (up to a plane isotopy). The reconstruction of $G$ is a reversed process of the subdividing and contracting (used to obtain $G_1$ and $G_2$): deleted loop edges are first glued back to $G_1$ and $G_2$, the contracted edges (if any) between two blue vertices $v_1^b$ and $v_2^b$ are expanded back, then the blue vertices in $G_1$ and $G_2$ resulting from the contraction of $\gamma^\p$ are expanded back to a closed curve equivalent to $\gamma^\p$ and the edges cut by $\gamma^\p$ are glued back together in the last step. This reversed process is made possible since all the information needed is stored in the gluing instruction.
We summarize this in the following lemma.

\begin{lemma}
\label{reconstructionlemma}
Let $G$ be a BRT-graph and let $G_1$ and $G_2$ be the two BRT-graphs obtained by a vertex-cut subdivision or an edge-cut subdivision of $G$.  Then the planar embeddings $G_1$ and $G_2$ induced from $G$ together with the gluing instruction that arise from this subdivision allow a reconstruction of a graph that is plane isotopic to the original graph $G$.
\end{lemma}

\section{Balanced Subdivisions of BRT-graphs and Knot Diagrams}\label{sec5}

\smallskip
In the last section, it was shown that subdividing a BRT-graph $G$ by a vertex-cut or by an edge-cut (based on a normal cut-cycle) results in two BRT-graphs $G_1$ and $G_2$. The main task of this section is to show that it is possible to subdivide a BRT-graph $G$ such that the sizes of $G_1$ and $G_2$ are balanced. Recall from the definition of the subdivision that a red vertex of $G$ remains a red vertex in one of $G_1$ and $G_2$ (but not both) and no new red vertices are created in the process. That is, if $G_1$ and $G_2$ are the graphs obtained from $G$ by subdivision with $V^R_1$ and $V^R_2$ being the sets of red vertices respectively, then $V^R$ is the disjoint union of $V^R_1$ and $V^R_2$. For a BRT-graph $G$, let us define its {\em standard weight} $W_s(G)$ as the number of its red vertices, i.e., $W_s(G) = |V^R|$. If $G$ is subdivided into $G_1$ and $G_2$ by a vertex-cut or an edge-cut, then $W_s(G) = W_s(G_1)+W_s(G_2)$. In this section non-standard weight systems are used which are denoted by a lower case $w$-function.

\begin{definition}
Let $G$ be a BRT-graph and let $c>0$ be a constant independent of $G$. A subdivision of $G$ into $G_1$ and $G_2$ (by either a vertex-cut or an edge-cut) is {\em balanced} if $\min\{W_s(G_1),W_s(G_2)\} \geq W_s(G)/6$ and in the case that the subdivision is an edge-cut subdivision, the length of the normal cycle (i.e., the number of red vertices in the cycle) used for the edge-cut is at most $c \sqrt{W_s(G)}$.
\end{definition}

\begin{figure}[htbp]
\begin{center}
\includegraphics[scale=0.6]{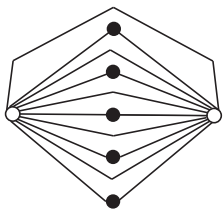}
\caption{ A BRT-graph with high blue/red vertex ratio. The red vertices are white and the blue vertices are black.}
\label{bluevertexbound}
\end{center}
\end{figure}

In a general BRT-graph, the number of blue vertices may not be bounded above by a function of the number of red vertices as shown in Figure \ref{bluevertexbound}. However, if the degrees of the red vertices of a BRT-graph $G$ are bounded above by a constant $g\ge 4$, then the number of blue vertices in $G$ is related to the number of red vertices in $G$ as shown in Lemma \ref{vbbound}.

\begin{lemma}
\label{vbbound}
Let $G$ be a BRT-graph with $|V(G)|\ge 4$ and $g$ be an upper bound of the degrees of the red vertices,
then $|V^B|\le (g/4) |V^R|$. Consequently, $|V(G)|\le (1+g/4)|V^R|$.
\end{lemma}

\begin{proof}
It suffices to prove the inequality for a BRT-component $G(V^*_M)$ of $G$. Deleting all but one edge from each set of multiple edges connecting the same two vertices in $G(V^*_M)$ results in a simple graph $H$ with the same number of blue and red vertices as that of $G(V^*_M)$.
(Here we think of multiple edges as edges that create an empty digon and not just edges that have the same end vertices. For example, the graph in Figure \ref{bluevertexbound} does not have a multiple edge.)
Now $H$ has only triangular faces. Let $n_r$ ($n_b$) be the number of red (blue)
vertices in $H$ and let $f$ be the number of faces in $H$.
The boundary of each triangular face contains at least 2 red vertices and any red vertex can be on the boundaries of at most $g$ different faces. Thus $n_r$ is bounded below by $n_r\ge 2f/g$. On the other hand, each blue vertex is on the boundaries of at least two faces and the boundary of each face contains at most one blue vertex (Lemma \ref{L3.3}). So $n_b\le f/2\le g n_r/4$.
\end{proof}

\medskip
\begin{lemma}\label{balancelemma}
Let $g\ge 4$ be a given constant. There exists a constant $W_0>3$ such that for any BRT-graph $G$ with $W_s(G)>W_0$ and the maximum degree of the red vertices in $G$ being $\le g$, there exists a balanced subdivision of $G$.
\end{lemma}

\begin{proof}
First consider the case that $G$ contains a BRT-component $G(V^*_M)$ such that $W_s(G(V^*_M))\ge W_s(G)/2$.
(Note that it is possible that $G=G(V^*_M)$.)  A non-standard weight system $w$ is assigned to $G(V^*_M)$ as follows. Let $m=|V^R|$ be the number of red vertices in $G$. Each red vertex in $G(V^*_M)$ is assigned weight $1/m$. All blue vertices in $G(V^*_M)$ are assigned weight zero.
Each face $f$ of $G(V^*_M)$ is assigned a weight $w(f)=r_f/m$, where $r_f$ is the number of red vertices of $G\setminus G(V^*_M)$ that are contained in $f$.

The total weight is equal to 1 since every red vertex of $G$ is either in $M$ or is contained in a face of $G(V^*_M)$. Since $W_s(M)=W_s(G(V^*_M))\ge W_s(G)/2$, no face of $G(V^*_M)$ has weight larger than $1/2$. Under this non-standard weight assignment, Theorem \ref{cyclecut} implies that there exists a $\frac{2}{3}$-cycle cut
that divides $G(V^*_M)$ (hence $G$) into two subgraphs.
Moreover the length of the cycle $\gamma$ used is at most $2\sqrt{2n}$ where $n=|V(G(V^*_M))|\le (1+g/4)m$ by Lemma \ref{vbbound}.
If $\gamma$ is not normal, then it can be modified into a normal cut-cycle in $G(V^*_M)$ by Lemma \ref{cutcycleexistance}. The normal cut-cycle $\gamma_1$ so obtained is shorter than $\gamma$, and each of the two subgraphs separated by it has at most $(2/3)m+2\sqrt{2n}$ red vertices. Now choose $W_0>0$ to be a constant large enough so that $(2/3)m+2\sqrt{2n}\le (2/3)m+2\sqrt{2(1+g/4)m}<(5/6)m$ holds for every $m>W_0$.

\medskip
Next consider the case that every BRT-component $G(V^*_M)$ has a standard weight $W_s(G(V^*_M))< W_s(G)/2=m/2$ (where $m=|V^R|$).
Let $T_G$ be the tree defined in Section \ref{sec3} (before Lemma \ref{CycleinBRT}).
A red vertex $v_M$ of $T_G$ that corresponds to a red component $M$ is assigned the weight $w(v_M)=W_s(M)/m$.
All blue vertices of $T_G$ are assigned weight zero. Notice that under this weight assignment, the total weight
is 1. Thus by Theorem \ref{generalcut} there exists a cut-vertex $v$ in $T_G$ such that each connected component
in $T_G\setminus \{v\}$ has a total weight less than or equal to $2/3$. If $v$ is a blue vertex in $T_G$, then the
blue vertex $u$ in $G$ corresponding to $v$ is a cut vertex and can apparently be used to  obtain a balanced vertex-cut subdivision of $G$.
(A vertex-cut using the cut vertex $u$ obtained in this manner in facts leads to connected components each of which has a weight of $2m/3$ or less.)
On the other hand, if $v$ is a red vertex then it corresponds to a red component $M$ of $G(V^R)$. Assign $G(V^*_M)$ the non-standard weight
system $w_1$ as before: each red vertex in $M$ is assigned the weight $1/m$, each blue vertex in $G(V^*_M)$ is assigned weight zero, and each face $f$ of $G(V^*_M)$ is assigned the weight $w(f)=r_f/m$, where $r_f$ is the number of red vertices of $G\setminus G(V^*_M)$ that are contained in $f$.
Again the total weight is 1 since every red vertex of $G$ is either in $M$ or is contained in a face of
$G(V^*_M)$. No face $f$ in $G(V^*_M)$ has a weight $w(f)> 2/3$ since otherwise deleting the red vertex $v$
in $T_G$ corresponding to $M$ results in a connected component in $T_G$ with weight $>2/3$, contradicting the given
property of $v$. Thus by Theorem \ref{cyclecut} there exists a cycle $\gamma$ in $G(V^*_M)$ that yields a
$\frac{2}{3}$-cycle cut of $G$. Again modify $\gamma$ as before to obtain a normal cycle $\gamma_2$
and use $\gamma_2$ to obtain an edge-cut subdivision of $G$. The only difference is that this time $\gamma_2$ causes a smaller bound on $W_0$ since the weight of $G(V^*_M)$ is less than $m/2$,  so the total weight of each
of the two graphs obtained by the edge-cut using $\gamma_2$ is bounded above by $(2/3)m+2\sqrt{2n}<(5/6)m$ where $n\le (1+g/4)m/2$ is the number of vertices in $G(V^*_M)$.
\end{proof}

The definition of a (balanced) vertex-cut may allow many different choices for $G_1$ and $G_2$ by choosing different unions of BRT-components. In order to allow a successful reconstruction, constraints are imposed on the selection of the BRT-components for $G_i$ for a balanced vertex-cut. Lemma \ref{diskcomponent} specifies these constraints and asserts that they can always be met.

\begin{figure}[htbp]
\begin{center}
\includegraphics[scale=0.8]{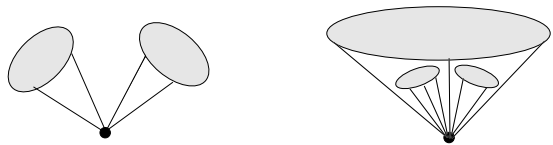}
\caption{Left: A separable disk-component; Right: An inseparable disk component containing two maximal disk-components. A gray area indicates BRT components as in Figure \ref{vertexcut}.}
\label{maxdisks}
\end{center}
\end{figure}

\medskip
Assume that $G$ is a BRT-graph and $v$ is a blue vertex in $G$ that can be used for a balanced vertex-cut. Let $\alpha$ be an arc that starts and ends at $v$ and is otherwise disjoint from $G$. $\alpha$ separates $G$ into two subgraphs $G_1$ and $G_2$ both containing $v$ and that are unions of complete BRT-components. We call $G_1$ (and $G_2$) a {\it disk-component} of $G$. It is possible that one of the two graphs contains only the vertex $v$. However, if both $G_1$ and $G_2$ contain at least one vertex other than $v$, then $G_1$ and $G_2$ are called {\it proper disk-components} of $G$. A disk-component $H$ is {\it separable} if there exist two proper disk-components $H_1$ and $H_2$ such that $H=H_1\cup H_2$ and $H_1$ lies in the outer face of $H_2$. $H$ is called {\it inseparable} if it is not separable. $H^\p$ is a {\it maximal disk-component} of $H$ if $H^\p$ is a proper disk-component contained in $H$ and if for every disk-component $D$ in $H$ that contains $H^\p$ either $D=H^\p$ or $D=H$. See Figure \ref{maxdisks}.

\medskip

\begin{lemma}\label{diskcomponent}
Let $G$ be a BRT-graph that admits a balanced vertex cut using the cut-vertex $u$ (in $G$), where $u$ is a cut vertex obtained as in the proof of Lemma \ref{balancelemma}. Let $G_1$ and $G_2$ be the two subgraphs obtained by the vertex-cut.  Then one of the subgraphs, say $G_1$, can be chosen as one of the following:

(i) $G_1$ is a disk-component of $G$.

(ii) $G_1$ consists of a union of several maximal disk-components of a single inseparable disk-component of $G$.
\end{lemma}

Notice that Figure \ref{vertexcut} is an example of case (ii).

\begin{proof}
If there exists a disk-component $H$ of $G$ with $W_s(G)/6 \le W_s(H)\le 5W_s(G)/6$ then let $G_1$ be $H$ (case (i)).
Now assume that no disk-component $H$ exists in $G$ such that $W_s(G)/6 \le W_s(H)\le 5W_s(G)/6$.
Assume there exists a disk-component $H$ of $G$ such that $W_s(H)> 5W_s(G)/6$. Moreover assume that among all disk-components $H$ with $W_s(H)> 5W_s(G)/6$, $H$ is the smallest one.

\medskip
Claim 1: $H$ must be inseparable. Otherwise, $H=H_1\cup H_2$ for two disjoint proper disk-components $H_1$ and $H_2$. $H_1$ or $H_2$ must have weight less than or equal to $5W_s(G)/6$ since $H$ is the smallest disk-component with weight more than $5W_s(G)/6$. But then it must be true that  $W_s(G)/6 \le W_s(H_1)\le 5W_s(G)/6$ or $W_s(G)/6 \le W_s(H_2)\le 5W_s(G)/6$ since $W_s(H_1)+W_s(H_2)=W_s(H)>5W_s(G)/6$. This is a contradiction since we assumed that there are no disk-components with weight between $W_s(G)/6$ and $5W_s(G)/6$.

\medskip
Claim 2: $H$ must contain at least one proper disk-component $H^\p$ with $W_s(H^\p)<W_s(H)$. If this is not the case, deleting $u$ from $G$ results in a connected component of weight more than $5W_s(G)/6$, contradicting the fact that $u$ is a cut-vertex for a balanced vertex cut. Remember that the cut-vertex $u$ obtained in the the proof of Lemma \ref{balancelemma} leads to connected components with weights $\le 2/3 W_s(G)$.

\medskip
It follows that $H$ contains maximal proper disk-components.
Let $H_1, H_2,\ldots,H_k$ be the maximal proper disk-components of $H$, then for each $i$, $W_s(H_i)<W_s(G)/6$ by our assumptions.
Let $W_s = \sum_i W_s(H_i)$. We must have $W_s \ge W_s(H)/6$, otherwise the graph $H\setminus (\cup_i H_i)$ has weight $>5W_s(G)/6-W_s(G)/6 =2W_s(G)/3$ and remains connected after $u$ is deleted, contradicting the fact that $u$ is a cut-vertex for a $2/3$-balanced vertex cut.
Thus $G_1$ can be chosen to be the union of some or all of the $H_i$'s (case (ii)).
The last case we need to consider is that all disk-components $H$ of $G$ satisfy the condition $W_s(H)< W_s(G)/6$. However this is impossible since $G$ is a disk-component of itself.
\end{proof}

\medskip
In order to apply the divide-and-conquer technique, it is necessary for us to use repeated balanced subdivisions to a BRT-graph $G$.

\begin{definition}\label{recursivesubdivision}
A BRT-graph $G$ (and the resulting BRT-subgraphs) can be divided recursively using balanced subdivisions. When the standard weight of a BRT-graph obtained in this repeated subdivision process falls below a pre-determined threshold $W_0$, the subdivision process stops on this BRT-graph and it is called a {\em terminal BRT-graph}. The subdivision process has to terminate at the point when all the resulting BRT-graphs are terminal BRT-graphs. The balanced subdivisions used to reach this stage are called a {\em balanced recursive subdivision sequence} of $G$.
\end{definition}

\begin{figure}[htbp]
\begin{center}
\includegraphics[scale=1.0]{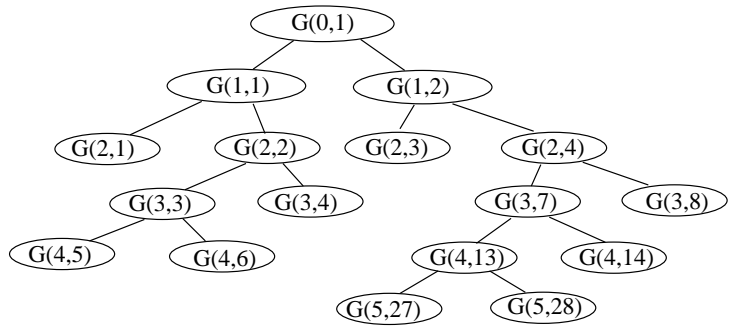}
\caption{The tree structure of BRT-graphs obtained from a balanced recursive
subdivision sequence of $G$.}\label{subdivisiontree}
\end{center}
\end{figure}

\medskip
To keep track of the BRT-graphs obtained when a balanced recursive subdivision sequence is applied to a BRT-graph $G$,  the following notations are adopted. $G(0,1)=G$. When the first subdivision is applied, the two resulting BRT-graphs are denoted by $G(1,1)$ and $G(1,2)$, the gluing instruction of the subdivision process is denoted by $N(0,1)$, and  the two new blue vertices created are denoted by $v^b(1,1)$ and $v^b(1,2)$. The two BRT-graphs obtained from subdividing $G(1,1)$ are denoted by $G(2,1)$ and $G(2,2)$ and the two BRT-graphs obtained from subdividing $G(1,2)$ are denoted by $G(2,3)$ and $G(2,4)$, and so on. In general, the BRT-graphs obtained from subdividing $G(i,j)$ (if $W_s(G(i,j))>W_0$) are denoted by $G(i+1,2j-1)$ and $G(i+1,2j)$, the newly created blue vertices are $v^b(i+1,2j-1)$ and $v^b(i+1, 2j)$, and the gluing instruction is $N(i,j)$. See Figure \ref{subdivisiontree} for an illustration of this relation. Notice that the lengths of the paths from the root ($G(0,1)$) of the tree to the leaves are not necessarily the same as shown in Figure \ref{subdivisiontree}, since some BRT-graphs may terminate earlier than others due to size differences. We say that a BRT-graph $H$ is {\em induced} from the plane graph $G$ if $H$ is one of the $G(i,j)$s described above. If $G(i_0,j)$ is a terminal BRT-graph where $i_0$ is largest among all other terminal BRT-graphs induced from $G(0,1)$ (from the same recursive subdivision sequence), then $i_0$ is called the {\em depth} of the corresponding recursive subdivision sequence.

\medskip
Lemmas \ref{splittinglemma} and \ref{balancelemma} lead to the following theorem.

\begin{theorem}\label{Trecexist}
There exists a balanced recursive subdivision sequence for each BRT-graph $G$. Furthermore, $W_s(G(i,j))\le |V^R(G)|(5/6)^{i}$ since
the subdivisions are balanced. It follows that there exists a constant $c_r>0$ ($c_r$ depends only on $W_0$ and the maximal degree $g$ of all the red vertices in $G$) such that the depth of any balanced recursive subdivision sequence of $G$ is bounded above by $c_r\ln(|V^R(G)|)$.
\end{theorem}

\begin{remark}\label{depthbound}
Let $i_0$ be the depth of a balanced recursive subdivision sequence of $G$ and $G(i,j)$ be one of the BRT-graphs induced from $G$ by this sequence. If we apply Theorem \ref{Trecexist} with $G(i,j)$ playing the role of $G$ (as the starting graph in the subdivision sequence), then the depth $d$ of the subdivision sequence leading $G(i,j)$ to its terminal BRT-graphs is at most $i_0-i$ and we have $d\le c_r\ln(W_s(G(i,j)))$.
\end{remark}

To apply the recursive subdivision to a knot diagram, we start with a minimum knot diagram $D$ of the knot $\K$ so that the number of crossings in $D$ is equal to $n=Cr(\K)$. Ignoring the over/under information of $D$ at its crossings, we treat $D$ as a 4-regular plane graph. In general, $D$ is not a BRT-graph since $D$ may contain faces of arbitrarily large size. Thus the previously established results cannot be applied directly to $D$. To remedy this problem, artificial edges are added to $D$ so that the resulting graph is a BRT-graph. These added edges may simply be removed from the embedding of the modified graph at the end of the process. The following lemma asserts that $D$ can be modified into a BRT-graph in such a way that the maximum degree of its vertices is bounded by a constant.

\begin{lemma}\label{deg12}
Let $D$ be a minimum projection of $\K$. If $D$ is treated as a plane graph so that crossings of $D$ are treated as vertices and strands connecting crossings are treated as edges, then by simply adding some new edges to $D$, $D$ can be modified into a plane graph $G$ such that $G$ is triangulated and the maximum degree of the vertices of $G$ is bounded above by 12.
\end{lemma}

\begin{figure}[htbp]
\begin{center}
\includegraphics[scale=0.4]{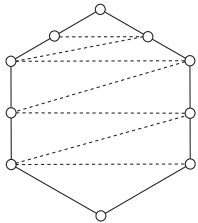}
\caption{ Triangulation of a face of a 4-regular plane graph by adding new edges. The edges added in the triangulation are dashed.}\label{facetriangulation}
\end{center}
\end{figure}

\begin{proof}
Each face $F$ of $D$ can be triangulated in a way as shown in Figure \ref{facetriangulation}. In doing so, at most two edges are added to a vertex of $F$. Since each vertex in $D$ belongs to
at most 4 faces this results in at most 8 new edges being added to each vertex. Hence the maximum degree of the resulting graph is bounded above by 12.
\end{proof}

Let $G$ be the triangulated graph obtained from $D$ as described in the proof of Lemma \ref{deg12}. At this stage, every vertex in $G$ is considered to be a red vertex. Since $G$ contains no blue vertices and no loop edges, it admits a proper BR-partition. Furthermore, it contains only one red component (namely itself) and this component is triangulated. Thus by definition, $G$ is a BRT-graph (without blue vertices). By Theorem \ref{Trecexist}, there exists a balanced recursive subdivision sequence for $G$. Since the subdivision operations do not increase the degree of a red vertex, the maximum degree of red vertices in such a graph $H$ is still bounded above by 12,  see Lemma \ref{deg12}.

\medskip
\section{Standard 3D-embeddings and Grid-like Embeddings of BRT-graphs}\label{sec6}

\medskip
In this section we introduce two special kinds of embeddings: the standard 3D-embedding and the grid-like embedding. The purpose of introducing the standard 3D-embeddings of BRT-graphs is to use these embeddings as benchmarks to verify that the topology of a graph is preserved when it is reconstructed from its two induced BRT-graphs. On the other hand, the purpose of introducing the grid-like embeddings is to simplify the reconstruction process: if two induced BRT-graphs are grid-like, then their grid-like structure will allow us to reconnect them in a way to preserve this grid-like structure so this reconnected graph can be used again in the next round of the reconstruction process of $G$. Furthermore, a grid-like embedding is almost on the lattice and in the last step when $G=G(0,1)$ itself is reconstructed (from its two immediate induced BRT-graphs $G(1,1)$ and $G(1,2)$), it will be easily modified into a lattice embedding.

\medskip
\subsection{Standard 3D-embeddings.}

In the following we assume that all BRT-graphs $G(i,j)$ involved are induced from a plane graph $G$ and that the maximum degree of red vertices is bounded above by 12.  Below we are introducing some terminology that we will use for the graphs $G(i,j)$ and their blue vertices throughout the next sections.

\medskip
\noindent
\textbf{Rectangles.}
Since $H$ = $G(i,j)$ is a plane graph drawn in the plane $z=0$, it can be embedded in the interior of a rectangle $R$ in the plane $z=0$. That is, there exists a plane isotopy $\Psi:\ \R^2\times [0,1]\longrightarrow \R^2$ such that $\Psi(x,0)=id$ and $H_1=\Psi(H,1)\subset (R\setminus \partial R)$. We will assume that all graphs $H=G(i,j)$ are contained in such a rectangle. Keep in mind that our the sides of our rectangles are parallel to either the $x$-axis or the $y$-axis.

\begin{figure}[htbp]
\begin{center}
\includegraphics[scale=0.8]{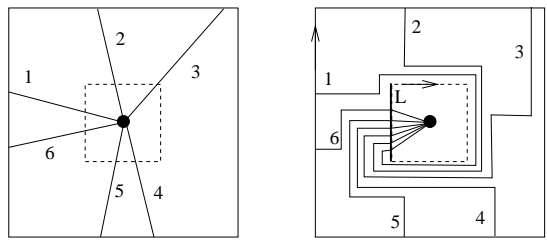}
\caption{Left: A blue square; Right: A blue square with edges inside it re-routed through a side $L$ (marked by the thickened line segment). The edge with label 1 is the first edge in the linear order assigned at the blue vertex. }
\label{bluesq_route}
\end{center}
\end{figure}

\medskip
\noindent
\textbf{Blue squares.} For each blue vertex $v$ of $H$ create a small square $S^b_v$ in the plane $z=0$ with side length $3\ell$ for some fixed small positive number $\ell >0$ such that $v$ is at the center of the square and $S^b_v\subset (R\setminus \partial R)$. $S^b_v$ is called a {\em blue square}. See Figure \ref{bluesq_route}.

\medskip
Without loss of generality we can assume that the boundary of $S^b_v$ intersects each edge leading out of $v$ exactly once transversely. Notice that  we can isotope the graph locally so that the edges within $S_v^b$ are single line segments as shown on the left of Figure \ref{bluesq_route}.

Let $R_v$ be the square with side length $\ell$ and center $v$ and $L$ be a side of $R_v$. The purpose of $R_v$ and $L$ is for us to use a local (VNP-)isotopy to re-route the edges connected to $v$ in such a way that they all enter $R_v$ from $L$. Moreover within $S_v^b$ the edges use only segments parallel to the $x$- and $y$-axis.  The right side of Figure \ref{bluesq_route} then shows an example of how to re-route these edges within $S^b_v\setminus R_v$ to achieve the desired result. As the figure shows, this can be done for each edge in $S_v^b$ involved with at most six right angle turns in the $xy$-plane. Let $e_v$ be the vector (parallel to either the $x$- or the $y$-axis) that points perpendicular from $L$ to $v$ and we call $e_v$ an {\em extension vector}.
\medskip

\noindent
\textbf{Orientation assignment of $\partial S^b_v$.}
Without loss of generality we can assume that the boundary of $S^b_v$ intersects each edge leading out of $v$ exactly once transversely. Recall that when a new blue vertex is created by an edge-cut subdivision, we oriented the components of $\partial N$ counterclockwise and used this orientation to define the cyclic order of the edges intersecting $\partial N$. If a new blue vertex $v=v^b(i,j)$ is in the BRT-graph $G(i,j)$ obtained by using the part of the original graph outside of $\gamma^\p$, then one may treat $\partial S^b_v$ as a deformation (contraction) of the outer component of $\partial N$. In this case we give $\partial S^b_v$ a counterclockwise orientation. On the other hand, if $v=v^b(i,j)$ is in the BRT-graph obtained using the part of the original graph inside of $\gamma^\p$, then  $\partial S^b_v$ should be treated as a deformation of the inner component of $\partial N$, where one would have to flip the inner component of $\partial N$ to realize the resulting graph on the plane without edge crossings. Thus in this case we will assign $\partial S^b_v$ a clockwise orientation. Finally, in the case that $v$ is created by a vertex-cut subdivision, $\partial S^b_v$ is always assigned the counterclockwise orientation.

\medskip
\noindent
\textbf{The orientation vector.} The order of the intersection points on $L$ is inherited from the order of intersection points on the boundary of $S_v^b$ induced by the orientation of $S_v^b$. To be more precise, we can choose any edge on the boundary of $S_v^b$ and using any path to connect it to $L$.  After that we can choose a second edge to go on either side of the first edge along $L$. After that the order of all other intersection points on $L$ is determined. In the example of Figure \ref{bluesq_route} we chose $L$ to be on the left side of $R_v$. Once we fix the edge with label 1 anywhere on $L$ then there are only two choices for the other edges to follow: we can obtain edge order $\{1,2,3,4,5,6\}$ ascending along $L$, or an edge order of $\{2,3,4,5,6,1\}$ ascending along $L$.  We can think of both of these as the same orientation along $L$ with the difference that one starts with the edge labeled 1 and the other with the edge labeled 2. Thus the counterclockwise or clockwise cyclic order of the edges around $S^b_v$ introduces a unique direction on $L$. We call the vector given by this orientation the {\em orientation vector} $o_v$. Furthermore, we will choose the first edge to intersect $L$ so that the linear order so obtained on $L$ matches the linear order at $v_b$, i.e. the linear order obtained from the counterclockwise orientation on $\partial N$ together with the choice of the path $\beta$ in $N$ (edge-cut) or the point $\beta$ on $\partial N$ (vertex-cut) as given by the gluing instructions.

\medskip
Notice that we use the same name and symbol for the small circular arrow around blue vertices in a neighborhood $N$ that contains the glueing instructions, see Definition \ref{glueinginstructions}.  We also use the word orientation vector of the orientation on $\partial N$. The reason is that these orientation vectors in $N$ or on $\partial N$ directly induce the orientation vector along the segments $L$. The vectors are equivalent (they tell us the linear order) and this justifies the identical names.

\medskip
It follows that the extension vector $e_v$ can be obtained from the orientation vector $o_v$ by a 90 degree clockwise or counterclockwise rotation. In the example of Figure \ref{bluesq_route} we need to rotate the the orientation vector $o_v$ 90 degree clockwise to obtain the extension vector $e_v$. This tells us that in this case $v$ arose in an edge-cut and belongs to $G_1$ (the graph obtained from the inside of $\gamma^\p$). If we need to rotate the orientation vector $o_v$ 90 degree counter clockwise to obtain the extension vector $e_v$ then the vertex $v$ arose in a vertex-cut or $v$ arose in an edge-cut and belongs to $G_2$ (the graph obtained from the outside of $\gamma^\p$).

\medskip
Finally, if in the formation of a blue vertex $v$ an edge between a newly temporarily created blue vertex and an older existing blue vertex is contracted, then the above orientation determination still applies: The linear order of the intersection points on $\partial N$ does include the edges from the existing blue vertex and these existing orders do not change when the blue vertices are merged.

\begin{figure}[htbp]
\begin{center}
\includegraphics[scale=0.8]{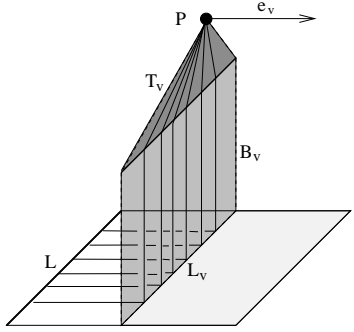}
\caption{An illustration of blue band and blue triangle: the construction only occurs within a rectangular box of height $t$ with $R_v$ as its base. The unused area in $R_v$ is lightly shaded.}
\label{bluetriangle}
\end{center}
\end{figure}

\medskip
\noindent
\textbf{Blue bands and blue triangles.} Let $L_v$ be the cross-section of $R_v$ that contains the blue vertex $v$ and is parallel to $L$ and let $L^\p_v$ be the line segment obtained by moving $L_v$ up to the plane $z=t_1>0$ for some positive $t_1$. Let $\delta>0$ be a small positive real number and let $P$ be the point directly over $v$ in the plane $z=t$ where $t=t_1+\delta$. The rectangle formed defined by $L_v$ and $L^\p_v$ is called a {\em blue band} and is denoted by $B_v$. The (vertical) triangle formed by $L^\p_v$ and $P$ is called a {\em blue triangle} and is denoted by $T_v$. We will now redraw the graph locally within the boundary of $R_v\times \mathbb{R}^+$ as shown in Figure \ref{bluetriangle}. Under this redrawing, the point $P$ becomes the blue vertex $v$ and each edge from the boundary of $S_v^b$ to $v_b$ is replaced a path consisting at most 9 straight line segments: at most 7 in the $xy$-plane, one vertical (from $L_v$ to $B_v$) and one slant (from $B_v$ to $P$). This obviously does not change the topology of the graph, it simply creates a 3-D structure of the graph for us to work with. We place the extension vector $e_v$ at the point $P$ for future references.
See Figure \ref{bluetriangle} for an illustration of this. Note that half of $R_{v_b}$ with $L_{v_b}$ as a side but opposite to $L$ is not occupied by any edge under this construction.

Assume that the above process is applied to every blue vertex in the graph, then we arrive at a new graph $H_s$ that is bounded in the rectangular box $B=R\times [0,t]$. Finally, we require that the projections of the blue squares to the $x$-axis and to the $y$-axis do not overlap each other. This can be done since we can pre-determine the positions of the blue vertices in $R$ by Lemma \ref{fixedcurvelemma} and we can choose the side length of these squares arbitrarily small.

\medskip
\begin{definition}\label{3ddrawing}
The embedding $H_s$ obtained in the above process from the BRT-graph $H$ in the rectangular box $B$ is called a {\em standard 3D-embedding} .
\end{definition}

Notice that a standard 3D-embedding $H_s$ of $H$ is isotopic to the BRT-graph $H$ by a VNP-isotopy that is the identity outside the space $R\times [-\epsilon,t+\epsilon]$ for an arbitrarily small positive constant $\epsilon$. Let us first describe the part of the isotopy that involves the space $z>0$. Here at  each blue vertex $v$ we retract the blue band $B_v$ until the base of the blue triangle reaches the plane $z=0$. After that we fold the blue triangle rigidly into the half of $R_{v_b}$ that is not occupied by any edge. This give us an embedding that is entirely contained in the plane $z=0$. The rest of the VNP-isotopy involves only moves within the plane (which is a plane isotopy that is automatically VNP). In other words, the VNP-isotopy described here is just a trivial extension of a plane isotopy.

\medskip
From Lemma \ref{redrawGinR} we know that we may pick the locations of the deformed blue vertices on the top of the box $R\times [0,t]$ with almost total freedom. Therefore, we like to have standard 3D-embeddings of $H$ with the property defined in the following definition.

\medskip
\begin{definition}\label{predetermined3ddrawing}
Let $H$ be a BRT-graph induced from $G$ with blue vertices $\lbrace v_1,\ldots, v_k\rbrace$. Let $R$ be any given rectangle in $z=0$ whose sides are parallel to either the $x$- or $y$-axis. Let $Q_1$, $Q_2$, ..., $Q_k$ be any $k$ distinct points in the interior of $R$ and let $P_1$, $P_2$, ..., $P_k$ be the corresponding points on the plane $z=t$ directly above the points $Q_1$, $Q_2$, ..., $Q_k$. If there exists a standard drawing $H_s$ of $H$ in $R\times [0,t]$ such that $P_j$ is the blue vertex in $H_s$ corresponding to $v_j$ for $j=1$, $2$, ..., $k$, then $H_s$ is called pre-determined standard 3D-embedding of $H$ (with $P_1$, $P_2$, ..., $P_k$ being the pre-determined blue vertices).
\end{definition}

The following lemma asserts that it is indeed possible to create pre-determined standard 3D-embeddings of $H$.

\begin{lemma}\label{pre3dexist}
Let $H$ be a BRT-graph induced from $G$ with $k$ blue vertices and let $R$ be any given rectangle in $z=0$ whose sides are parallel to either the $x$- or $y$-axis. Let $Q_1$, $Q_2$, ..., $Q_k$ be any $k$ distinct points in the interior of $R$ and let $P_1$, $P_2$, ..., $P_k$ be the corresponding points on the plane $z=t$ directly above the points $Q_1$, $Q_2$, ..., $Q_k$. Then there exists a pre-determined standard 3D-embeddings of $H$ with $P_1$, $P_2$, ..., $P_k$ being the pre-determined blue vertices.
\end{lemma}

\begin{proof}
Let $H$ be a BRT-graph induced from $G$ with blue vertices $\lbrace v_1,\ldots, v_k\rbrace$. Let $R$ be any given rectangle in $z=0$ whose sides are parallel to the $x$- and $y$-axis. Let $Q_1$, $Q_2$, ..., $Q_k$ be any $k$ distinct points in the interior of $R$.  Then by Lemma \ref{redrawGinR} there exists a plane isotopy $\Psi$ such that $\Psi(G,1)$ is contained in $R$ and $\Psi(v_j,1)=Q_j$.  We can then obtain the desired pre-determined standard 3D-embedding
of $\Psi(G,1)$ with the $P_j$s being the blue vertices of the new graph by the previously outlined construction.
\end{proof}

\subsection{Grid-like Embeddings.}\label{sec7}
For the purpose of embedding the graph $G$ into the cubic lattice, the structure offered by a standard 3D-embedding is not enough. We need to use a structure that is almost like a lattice embedding for the graphs obtained in the subdivision process. A graph embedding with such a structure will be called a {\em grid-like} embedding. The detailed description of this embedding is given in this section.

\medskip
Assume that $H$ = $G(i,j)$ for some valid $i,j$ from the recursive subdivision process and that the degree of any red vertex in $G$ is at most 12.

\begin{definition}\label{gridlikedrawing}
We call an embedding $H_{gr}$ of a BRT-graph $H$ a {grid-like embedding} if it satisfies the following conditions:

\smallskip\noindent

(i) All red vertices of $H_{gr}$ are lattice points in the plane $z=0$. Moreover, each red vertex $v$ is contained in the interior of a lattice rectangle $S^R_v$ of dimensions $w \times l$ where $w, l\ge 3$  (called a {\em red square}) that does not contain any other vertices of $H_{gr}$. The edges of $H_{gr}$ connected to $v$ must pass through (different) lattice points on the boundary of $S^R_v$. Figure \ref{redsquare} shows this for a vertex of degree 12 in the smallest possible lattice rectangle (a  $3 \times 3$ square).

\begin{figure}[htbp]
\begin{center}
\includegraphics[scale=0.3]{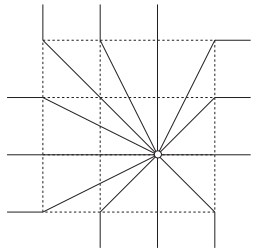}
\caption{ A red vertex of degree at most 12 can be realized in a (red) $3\times 3$ lattice square.}
\label{redsquare}
\end{center}
\end{figure}

\smallskip\noindent
(ii) $H_{gr}$ is contained in a rectangular box $B=R\times [0,t]$ for some integer $t>0$. All blue vertices of $H_{gr}$ are on the top of the box, i.e., in $R\times \{t\}$. Similar to the standard 3D-embedding, each blue vertex is the top vertex of a blue triangle that is on top of a blue band originated from a blue square. The difference here is that the blue square (from where the blue band originates) may be in a horizontal plane $z=h$ for some integer $0< h<t-\delta$ (it can also be in the plane $z=0$), here $\delta>0$ is a small positive constant. All the other requirements on the blue bands, triangles, and squares as given in the definition of a standard 3D-embedding must also be satisfied.

\smallskip\noindent
(iii) All edges outside the red squares and blue triangles are paths consisting of only line segments parallel to the $x$-, the $y$-, or the $z$-axis. Moreover all horizontal line segments must have integer $z$-coordinates.

\smallskip\noindent
(iv) All red-red edges (edges connecting two red vertices) are on the (cubic) lattice with the (possible) exception of the segments contained in red squares. The red-blue edges do not have to be on the lattice.

\smallskip\noindent
(v) $H_{gr}$ is isotopic to a standard 3D-embedding
of $H$ by a VNP-isotopy $\phi$ that is identity outside the box $R\times [-1,t-\delta]$. Furthermore the isotopy restricted to any red square must be a plane isotopy, that is throughout the isotopy a red square remains in the plane $z=0$. However it does not need to remain a lattice rectangle, the red square simply will play the role of a disk neighborhood of the red vertex.
\end{definition}

\medskip
The requirement that the VNP-isotopy is the identity outside the box $R\times [-1,t-\delta]$ enforces that the edge order of the blue vertices cannot be changed since the blue triangles do not move at all. The requirement on the red squares enforces that the edge order of the red vertices cannot be changed either. The reason for the extension $-1\le z< 0$ of the space the isotopy can use will become clear later.

\medskip

\section{Grid-like Embeddings of BRT-graphs Induced from a Knot Projection}\label{sec77}

\medskip
There are two approaches to obtain a grid-like embedding of a BRT-graph $H$. The first is a direct construction from the graph $H$ and the second is a reconstruction using two grid-like embeddings of the two BRT-graphs obtained from $H$ by a vertex-cut or an edge-cut subdivision. For the terminal BRT-graphs obtained from $G$ in the subdivision process, we will have to use the first approach to obtain their grid-like embeddings. Although we can use this first approach to get a grid-like embedding of $G$ itself as well, it will not achieve the desired efficiency in the embedding length. For that we will then need the second approach to assemble these grid-like embeddings of the terminal BRT-graphs into grid-like embeddings and ultimately obtain a grid-like embedding of $G$ (which will then be modified into a lattice embedding of the knot diagram).

\medskip

\subsection{Grid-like embedding via direct construction.} The following lemma assures that the first approach is always possible.

\begin{lemma}\label{pregridexist}
A grid-like embedding of a BRT-graph $H$ can be directly constructed from $H$.
\end{lemma}

\begin{proof}
Let $H_1$ be a standard 3D-embedding of $H$ that is guaranteed by Lemma \ref{pre3dexist}. We stretch $R$ in the $x$- and $y$-direction by inserting additional lattice lines. First we add lattice lines at the locations of the red vertices to put these on lattice. Next we add enough lattice lines to create required minimum size red rectangle.

\medskip
By induction on the number of edges in $H_1$, we can prove that all edges, except the parts contained in the red rectangles $S_i$'s can be straightened by a VNP-isotopy so that they consist of only line segments parallel to a coordinate axis.  We now stretch $R$ in the $x$- and $y$-direction in a recursive manner. Each stretch keeps the line segments already on the square lattice on the lattice, but takes at least one line segment on a red-red edge path that is not on the square lattice to the lattice.

\medskip
We need to be careful to not disturb the existing blue squares and red rectangles by these stretches. If a blue square is intersected by a line $t$ that contains a segment of a red-red edge then we first move the blue square slightly to ensure that $t$ cuts edges in the blue square only transversely. Next the lattice line is inserted at $t$ which stretches one side of the blue square to a length of $t+3\ell$. The blue square is reconstructed inside this rectangle with side length $3\ell$ without intersecting $t$. See Figure \ref{bluesquareexpand} for an illustration of this.  If a red rectangle is intersected by a line $t$, adding the lattice line results in a red rectangle which still satisfies the required size restrictions and no additional steps must be taken.

\medskip
Furthermore, the red rectangles remain disjoint and the blue triangles still share no common $x$- or $y$-coordinates after the stretches. It is easy to see that this is always possible. The resulting graph is denoted by $H^{g}$.
This process does not change the structure of the blue bands, it only moved some blue squares with their corresponding blue bands and blue vertices by a rigid motion. Thus the graph $H^{g}$ satisfies all the conditions of a grid-like embedding.
\end{proof}

\medskip
\subsection{Grid-like embedding via re-connection.}\label{sec72}

\medskip
Using a recursive subdivision process of $G$ to construct a grid-like embedding of $G$ on the lattice requires a procedure to combine two grid-like embeddings of BRT-graphs (obtained from either a circular edge-cut or a vertex-cut subdivision) into a new grid-like embedding. More precisely, let $G_0=G(i,j)$, $G_1=G(i+1,2j-1)$ and $G_2=G(i+1,2j)$ be three BRT-graphs obtained in the subdivision process of $G$ by a vertex-cut subdivision or an edge-cut subdivision as defined in the paragraph after Definition \ref{recursivesubdivision}.
We assume that for $i=1,2$, $G_i$ has a grid-like embedding $G^g_i$ embedded in a box $B_i$, which has the rectangle $R_i$ as its base in $z=0$ and the height $t_i$.
In addition, we assume the following convention: if $G_0=G(i,j)$ was divided using a circular edge-cut, then $G_1$ refers to the BRT-graph which is derived from the graph inside $\gamma^\p$; if $G(i,j)$ was divided using a vertex-cut, then $G_1$ refers to the subgraph which contains the maximal disk-components mentioned in \ref{diskcomponent} (ii), or if the balanced vertex cut is chosen following \ref{diskcomponent} (i) it refers arbitrarily to any one of the two subgraphs.

\medskip
This section describes a procedure to obtain a grid-like embedding $G^g_0$ of $G(i,j)$ using only information from the given embeddings $G^g_1$ and $G^g_2$ and from the gluing instruction $N(i,j)$. During the description, we refer to the blue vertices $v^b(i+1, 2j-1)$ and $v^b(i+1, 2j)$ created by the subdivision as $v^b_1$ and $v^b_2$, respectively.

\medskip
We split the construction into 7 steps.  Steps (1) through (3) serve to prepare the graphs $G^g_1$ and $G^g_2$ for the connection process, while
steps (4) through (7) make the actual connection of the edges.

\begin{figure}[htbp]
\begin{center}
\includegraphics[scale=1.0]{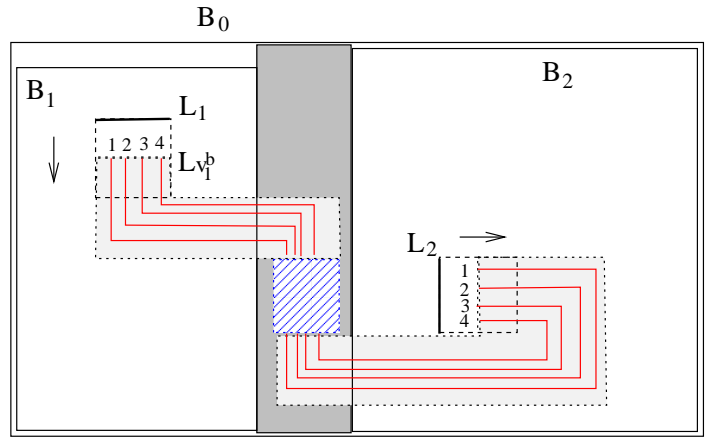}
\caption{The top view of the boxes $B_0$, $B_1$, $B_2$ and the connecting rectangle (shaded rectangle in the middle) in the case that the blue vertices are created by a circular edge-cut subdivision. Several typical paths from the top of a blue band to the other blue band are shown as well, with the labels of the edges to show how the corresponding edges should line up on the opposite sides of the connecting rectangle. Only the center parts of the blue squares are shown and they are greatly enlarged to reveal the details. }\label{bandex1}
\end{center}
\end{figure}

\medskip
(1) Align the boxes $B_i$ containing the $G^g_i$ properly in the $x$- and $y$-directions next to each other, with their base rectangles in the plane $z=0$. Without loss of generality we assume that the space between the two boxes is exactly one unit in the $x$-direction and that one side of the boxes coincides with the $x$-axis. Furthermore all the red squares and all the red-red edges outside the red rectangles are still on the cubic lattice. In the gap between the boxes we put a {\it connecting rectangle} at height $z=t_c = \max\lbrace t_1, t_2\rbrace$.  The connecting rectangle has dimensions $1\times y$ where $y = \max \{y_1,y_2\}$ and $y_i$ is the $y$-dimension of the box $B_i$ in their new location, see Figure \ref{bandex1}.

\medskip
(2) Create a (rectangular) box $B_0$ which includes both boxes $B_1$ and $B_2$ and the connecting rectangle and is of height $t_0 = t_c + 1$.

\medskip
(3) Extend the blue vertices in $G^g_1$ and $G^g_2$ other than the $v^b_i$ to the top of the box $B_0$ by extending their blue bands vertically by one unit and also lifting the blue triangles vertically by one unit.

\medskip
(4) Delete the blue triangles at the blue vertices $v^b_1$ and  $v^b_2$ and extend the corresponding blue bands in the $z$-direction to the plane $t_c$. Then extend the blue bands horizontally in the plane $z=t_c$ to the connecting rectangle. We refer to this horizontal extension of the blue band as the {\it extension band}. The two extension bands and the connecting rectangle is referred to as the {\it connecting strip}.
The extension band consists of rectangles in the plane $z=t_c$. It starts in the direction of the extension vector in the case of circular edge-cut and in the opposite direction of the extension vector in the case of vertex-cut. With at most two right angle turns within $S\times \{t_c\}$ (where $S$ is the blue square of the corresponding blue vertex), it can be made moving toward the connecting rectangle. With two right angle turns in the connecting rectangle and a suitable bandwidth change, it is then connected to the extension band coming from the other blue vertex, see Figure \ref{bandex1} (which is a case of the circular edge-cut subdivision). Since the turns only happen in the connecting rectangle and in the blue squares (at the $z=t_c$ level), the projection of the extension band into the $xy$-plane does not intersect any other blue square (hence itself will not intersect any other blue band in the rectangular boxes $B_i$). This is true because of the properties of the blue squares and blue bands.

\begin{figure}[htbp]
\begin{center}
\includegraphics[scale=0.7]{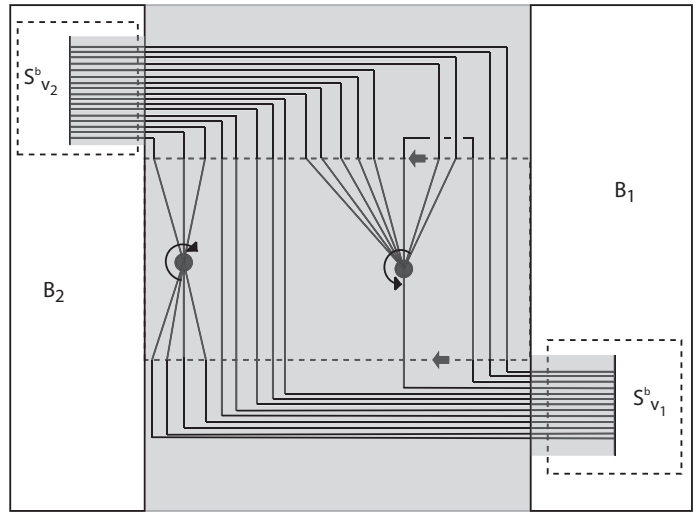}
\caption{How $G_1$ and $G_2$ are reconnected: The connecting rectangle is made much wider to show the details. In the middle of the connecting rectangle is $N^{\p\p}$, which is a deformation of the neighborhood of Figure \ref{unicycleneighborhood}. The curve that went under two edges represents a loop edge deleted when $G_1$ and $G_2$ were created.}\label{connectingstrip}
\end{center}
\end{figure}

\medskip
(5) This step applies when $G_1$ and $G_2$ are obtained after a circular edge-cut subdivision, so the extension band starts in the direction of the extension vector. Recall that we had assumed that in this case $G_1$ is the BRT-graph inside $\gamma^\p$. So $v^b_1$ is clockwise and $v^b_2$ is counterclockwise, see the paragraph marked as ``Orientation determination" before Definition \ref{3ddrawing}. Cut the neighborhood $N(i,j)$ open along $\beta$ and stretch it into a rectangle called $N^\p(i,j)$ as we did in Figure \ref{unicycleneighborhood}. Modify and re-scale $N^\p(i,j)$ so that (a) its side lengths are smaller than one third of the side lengths of either of the two blue squares and (b) each path connecting two opposite boundary points that correspond to the intersection points of the same edge with $\partial N(i,j)$ is just a single line segment (parallel to either the $x$- or the $y$-axis). This rectangle is then placed into the connecting rectangle and is denoted by $N^{\p\p}$. Since the linear order of the edges along the extension band is the same as the linear order of $v_1^b$ (and $v_2^b$) by the definition of grid-like embedding and the fact that $G_1$ and $G_2$ are grid-like embeddings. Thus edges with the same labels (edges that are to be connected) from each side (namely either from the $G_1$ or $G_2$ side) can be aligned perfectly with their counter parts of the edges of $N^{\p\p}$ as shown in Figure \ref{connectingstrip}.

\medskip
This connects all the edges in the blue bands and their extension bands arising from $G^g_1$ and $G^g_2$.  In particular, an edge from $G^g_1$ is connected to an edge $G^g_2$ only if both edges have the same label and no additional crossings are introduced.
By our construction, these edges consist  of only straight line segments. In particular any edge passing from the extension band of $G^g_1$ to the extension band of $G^g_2$ contains only two right angle turns on the connecting rectangle.

\begin{figure}[htbp]
\begin{center}
\includegraphics[scale=0.6]{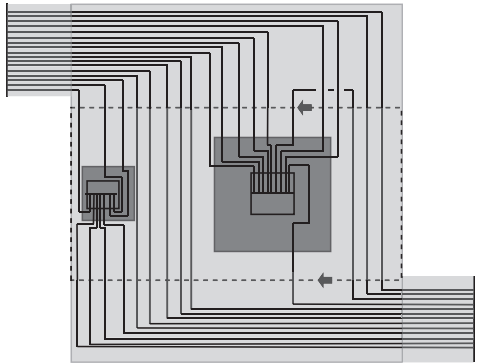}
\caption{The edges on the connecting rectangle of Figure \ref{connectingstrip} are grid-like and two new blue squares with blue bands and triangles have been created. The Figure is not to scale and the new blue squares appear to be much larger than they actually are. The linear orders at the new blue vertices are recovered from their corresponding gluing instructions. Only the top view is shown so the new blue vertices and blue triangles are not visible in the figure.}
\label{newbluesquares}
\end{center}
\end{figure}

\medskip
$N^\p(i,j)$ may contain one or more new blue vertices. For each such blue vertex $u$, extend each edge connected to $u$ from where it enters $N^{\p\p}$ to the centerline of $N^{\p\p}$, then extend it up by one unit. Create a small blue square in $z=t_c+1$ over $N^{\p\p}$ for this blue vertex. Then a blue band, and a blue triangle with the new blue vertex on top of the blue triangle in the plane $z=t_c+2$, following the same rules as before (for grid-like embeddings). See Figure \ref{newbluesquares}.

\medskip
Since the size of the new blue squares can be arbitrarily small and each blue vertex in $N^{\p\p}$ has certain free space to move (without crossing the straight paths that have been placed), the projections of the blue squares to the $x$- and the $y$-axes can be adjusted so as not to overlap with each other or with any other existing blue squares. For each blue square, the edges are combined into a vertical blue band as before and topped with a blue triangle at $z=t_c+1$. At this point all edges that are connected to $v^b_1$ or $v^b_2$ are accounted for.

\medskip
Note that we have not addressed the labels on $N(i,j)$ that correspond to loop edges that were deleted from $G_{1l}$ and $G_{2l}$ when $G_1$ and $G_2$ were created (although one such edge has been illustrated in Figure \ref{newbluesquares}). We will address this in Step (6).

\begin{figure}[htbp]
\begin{center}
\includegraphics[scale=0.8]{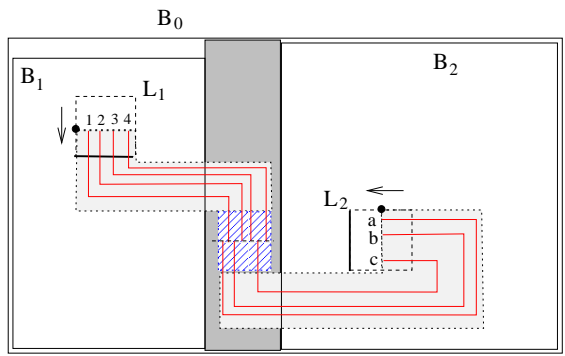}
\caption{Top view of an example for case 5$^\p$. In this case the original linear order of the edges at the blue vertex $v^b_0$ is $a1b23c4$. The solid dots can be thought of as the point $\beta$ used to define the linear order. A simplified 3D view of this is shown in Figure \ref{principalvertexcut}.}\label{bandex3}
\end{center}
\end{figure}

\begin{figure}[htbp]
\begin{center}
\includegraphics[scale=1.0]{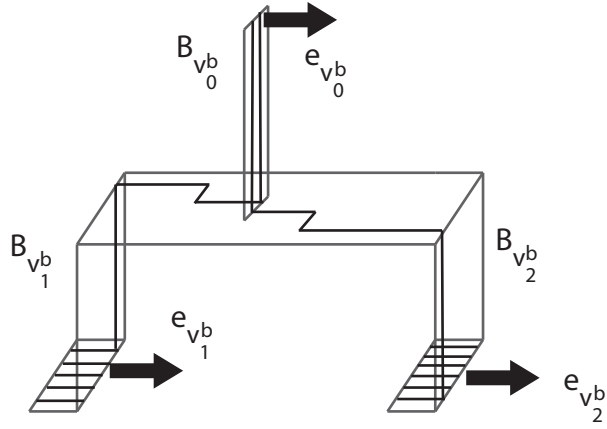}
\caption{ Shown are the three blue bands $B_{v^b_1}$, $B_{v^b_2}$, and $B_{v^b_0}$ (in gray) together with the three extension vectors $e_{v^b_1}$, $e_{v^b_2}$, and $e_{v^b_0}$ that arises in a vertex cut. The connection is only shown schematically without horizontal turns and the width of the bands is not to scale. Only two edges (black) are shown to illustrate the layout.}\label{principalvertexcut}
\end{center}
\end{figure}

\medskip
(5$^\p$) This step applies when $G_1$ and $G_2$ are obtained after a vertex-cut subdivision. For the vertex-cut, $v_1^g$ and $v_2^g$ have the same orientations in the plane. Since the extension band in this case starts in the direction opposite to the orientation vector, the linear orders of the edges along the band sides will again align correctly. See Figures \ref{bandex3} and \ref{principalvertexcut} for an illustration of this.

\medskip
The connecting rectangle contains $N(i,j)$ with a single blue vertex $v^b_0$ that was split into the two blue vertices  $v^b_1$ and $v^b_2$. As in the case of an edge-cut, a blue square is created in $z=t_c+1$ first, then a blue band, and a blue triangle with the new blue vertex on top of the blue triangle in $z=t_c+2$, satisfying all requirements of a grid-like embedding. This connects all the edges in the extension bands originated in the blue bands $B_{v^b_1}$ and $B_{v^b_2}$. Notice that the linear order of the edges at $v^b_0$ is restored when the edges from two sides meet the middle bar of the connecting rectangle.

%\begin{figure}[htbp]
%\begin{center}
%\includegraphics[scale=0.5]{vertexconnectingstrip.eps}
%\caption{ Shown is the connecting rectangle as it arises from the example in Figure \ref{vertexcut}. Shown are the two extension bands originating from $B_{v^b_1}$ and $B_{v^b_2}$, and the blue square $S_{v^b_0}$ together with $L_{v^b_0}$. The connection is only shown schematically. The width of the bands is not to scale. CHANGE FIGURE TO FORMAT OF FIG 12 and 13?}
%\label{vertexconnectingstrip}
%\end{center}
%\end{figure}

\medskip
(6) This step deals with the loop edges deleted after a circular edge-cut subdivision (such deletion can only happen in the case of a circular edge-cut subdivision). In the case that $G_1$ and $G_2$ are obtained after some loop edges are deleted, then $N^\p(i,j)$ contain labels not used to $G_1^g$ and $G_2^g$. The creation of temporary loop edges can happen in two ways.

\medskip
The first case is when a red-red edge $e$ (not on $\gamma$) was cut twice by $\gamma^\p$. Both red vertices incident to $e$ are contained in one of the graphs, say in $G_1$ (the case if the vertices are contained in $G_2$ is identical) and therefore in $G^g_1$. The middle arc of $e$ is contracted into a loop in $G^\p_2$ and is eventually deleted. This leads to 4 identical labels on $\partial N(i,j)$, two each on each component of $\partial N(i,j)$. Two labels in one boundary component of $\partial N(i,j)$ are accounted for by their corresponding edges on the extension band originated from $B^b_{v^b_1}$. The other two labels on the other boundary component of $\partial N(i,j)$ have no corresponding labels on $v^b_2$. The two edges coming up from $G^g_1$ just end on the opposite side of $N^{\p\p}$ (and there will be no other edges with the same label later on to connect them). However this is no problem since the gluing instruction tells us that these two edge ends must be connected to each other at this stage. Usually, this cannot be done in the plane $t_c$ without creating crossings. The connection is made using five edge segments, three of which are in the plane $t_c - 1$. Both edges are extended to the end of the connecting rectangle, see Figure \ref{newbluesquares} for an example of such an edge. A vertical segment of unit length is added at the end of each of the edges, connecting the plane at level $t_c$ with the plane at level $t_c -1$. A short segment parallel to the $y$-axis is added to the ends of both vertical segments and then one segment parallel to the $x$-axis connects the two end points. This construction builds a small `hook' which hangs below the connecting rectangle. Several such loops may have been removed during a subdivision step. For each of them a hook creates the correct connection between the edges without adding unwanted crossings. Loops may be nested, and the hooks can be nested too. A higher nesting level of the loops leads to slightly longer pairs of parallel segments which are parallel to the $y$-axis and a longer horizontal segment parallel to the $x$-axis.

\medskip
In the second case, a red-blue edge (not on $\gamma$) connected to the same blue vertex $v$ not on $\gamma$ are cut once by $\gamma^\p$. If there is just one such edge then this edge is contracted and does not lead to a loop. However if there is more than one, only one is contracted and the others result in loops. Let $e$ be one such edge. Then the red vertex connected to $e$ is contained in one of the graphs, say in $G_1$ (the case if the vertex is contained in $G_2$ is identical) and therefore in $G^g_1$. The part of $e$ connected to the blue vertex is contracted into a loop in $G^\p_{2l}$ and was deleted. This leads to 3 identical labels on $\partial N(i,j)$, one on the boundary component belonging to $G_1$ and the other two on the boundary component belonging to $G_2$. The single label in one boundary component of $\partial N(i,j)$ is accounted for by
by a corresponding label of edges on the extension band originating from $B^b_{v^b_1}$.  The two labels on the other boundary component of $\partial N(i,j)$ are not accounted for by corresponding labels on edges connected to $v^b_2$. We now construct a `small' hook exactly as in the first case.

\medskip
Notice that the loops in $G^\p_{1l}$ and $G^\p_{2l}$ are positioned at different ends of the neighborhood rectangle, that is,
the end that is closer to $B^b_1$ and $B^b_2$, respectively. The loops on each side nest perfectly but the loops on both sides combined may not exhibit a nesting behavior.
After this step the connecting rectangle contains edges accounted for all labels on $N(i,j)$ and in the same arrangement as is specified in $N(i,j)$.

\begin{figure}[htbp]
\begin{center}
\includegraphics[scale=0.7]{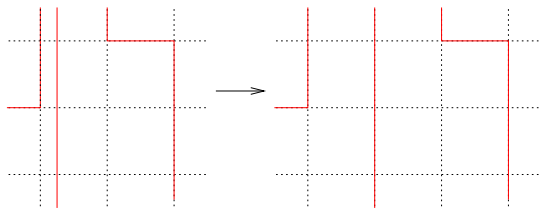}
\caption{Stretching the space to accommodate a newly created red-red edge on the lattice.}
\label{gridlineaddition}
\end{center}
\end{figure}

\medskip
(7) This last step only applies to the case when $G_1$ and $G_2$ are obtained after a circular edge-cut subdivision. The reason is in this case, we may have created red-red edges that are no longer on the lattice hence the reconnected graph described in the earlier steps is not grid-like yet. We will remedy this problem by adding new gridlines in the $x$- and $y$-direction to put newly formed red-red edges on lattice. For each line segment (on the re-connected red-red edge) parallel to the $y$-axis, a new $x$-gridline is added (which corresponds to a stretching isotopy). For each line segment parallel to the $x$-axis, a new $y$-gridline is added. See Figure \ref{gridlineaddition}.

\begin{figure}[htbp]
\begin{center}
\includegraphics[scale=0.7]{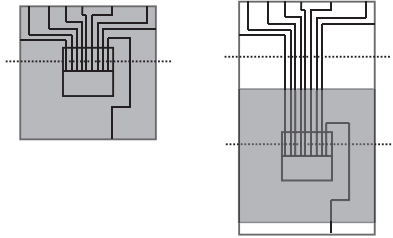}
\caption{On the left a blue square with a needed new gridline (dashed). On the right the blue square has been expand by one unit in the $y$-direction without creating any additional turns. The original center square $R_v$ has been translated in $y$-direction and a new blue square of the same size as the original has been created with $R_v$ at its center. The Figure is not to scale and the new blue squares appear to be larger than they are.}
\label{bluesquareexpand}
\end{center}
\end{figure}

\medskip
In this process we destroy the blue squares that formerly belong to the now vanished blue vertices $v^b_1$ and $v^b_2$. These blue squares are no longer needed. All other blue squares must be preserved.
If one of these new gridlines hits the projection of a blue square $S^b_v$ then we need to adjust the blue square as follows, see Figure \ref{bluesquareexpand}. First we slightly move the blue square to make sure that the new gridline hits edges in the blue square only transversely. Then we expand the blue square as is required by the insertion of the new grid line. The expanded $S^b_v$ now becomes a rectangle with one of its sides having a length of more than one unit. In this rectangle we put a copy of the original center square $R_v$ with the original small width $\ell$ by translation. (If we inserted an $x$- or $y$- gridline then we translate $R_v$ in $x$- or $y$-direction, respectively.) This can be connected up with exactly an many turns for the edges as before. Around this newly positioned square $R_v$ we reposition an new blue square $S^b_v$ with the original size. By default this new blue square has a projection that is disjoint in $x$- and $y$-coordinates from the projections of all the other blue squares.

\medskip
If one of these new gridlines hits a red rectangle then we simply stretch the red rectangle into a larger rectangle. This does not introduce any new turns and preserves all the required properties of a grid-like embedding.

This completes the description of how to combine the two grid-like embeddings of $G_1$ and $G_2$ into a new grid-like embedding of $G_0$.

\medskip
\section{The Verification of Topology Preservation}\label{sec73}

\medskip
While a grid-like embedding of a BRT-graph obtained as given in Lemma \ref{pregridexist} preserves its topology by its construction, it is far from obvious that the embedding obtained by reconnecting two grid-like embeddings as given in Section \ref{sec72} preserves the topology of the original BRT-graph from which the two grid-like embeddings are induced. We will prove that this is indeed the case in this section.

\medskip
\begin{lemma}
\label{extconstr}
Let  $G_1=G(i+1,2j-1)$ and $G_2=G(i+1,2j)$ be BRT-graphs obtained from the BRT-graph $G_0=G(i,j)$ by a subdivision. Then the grid-like embedding $G^g_0$ of $G_0$ as described in Section \ref{sec72}, confined in a rectangular box of the form $R^\p_0\times [0,t_0]$, is isotopic to $G(i,j)$ by a VNP-isotopy that is the identity outside a small neighborhood of the box  $B_0^\p=R^\p_0\times [-1,t_0-\delta]$. Therefore, $G^g_0$ is indeed a grid-like embedding of $G_0$ since it satisfies all other requirements of a grid-like embedding of $G_0$.
\end{lemma}

\medskip
Before stating the proof, let us recall that $G_1$ is either the interior graph in the case of a circular edge-cut subdivision, or the graph containing one disk-component or a union of some maximal disk-components in the case of a vertex-cut subdivision.

\begin{proof}
All the requirements for a grid-like embedding of $G(i,j)$ as specified in Definition \ref{gridlikedrawing} are already satisfied by the construction process outlined in Section \ref{sec72}, as one can check. Thus it suffices to show that the grid-like embedding obtained is isotopic to a standard 3D-embedding of $G_0$.

\medskip
In the last step in the construction process described in Section \ref{sec72}, the box $B_0$ (and the boxes $B_1$ and $B_2$) are stretched  to make room for new gridlines so that the newly created red-red edges can be put on the lattice. Let us call the stretched boxes $B_0^\p$, $B_1^\p$ and $B_2^\p$.

\medskip
Case 1: The case of a circular edge-cut subdivision.  The VNP-isotopy will be constructed by a sequence of isotopies using the following steps:

\medskip\noindent
(1) Remove some gridlines to get back to an almost grid-like embedding;

\smallskip\noindent
(2) Shrink $B_1^\p$ and $B_2^\p$ (with the graphs contained in them) back to $B_1$, $B_2$ and deform the graphs contained in them to 3D standard embeddings;

\smallskip\noindent
(3) Deform the hooks (representing the loop edges that were deleted during the subdivision) from the connecting rectangle into the plane $z=0$;

\smallskip\noindent
(4) Shrink $G_1$ and drag the shrank $G_1$ along the extension bands and the connecting rectangle and drop it into $G_2$ on $z=0$.

\smallskip\noindent
(5) Straighten out the blue squares and make the blue-bands vertical as required by a standard 3D-embedding of $G_0$ using Lemma \ref{straighten}.

\medskip
We now address each of these points in detail.

\medskip
(1)  Clearly such stretching isotopies are reversible. However the definition of a grid-like embedding does not allow us to use an isotopy that changes the entire box $B_0$. Thus we can only shrink the box $B_0^\p$ (together with $B_1^\p$ and $B_2^\p$) back to its original size under $z=t_0-2=t_c$. The first isotopy is defined by this shrinking isotopy $\phi^\p$ for $z\le t_0-2$, the identity for $z\le -1$ and $z\ge t_0-1-\delta$ where $\delta>0$ is the number chosen so that all blue triangles have bases on $z=t_0-1-\delta$. Let $f$ be a blue band at a blue vertex $v$. $f$ intersects the base of its corresponding blue triangle at $b$ and intersects the plane $z=t_0-2$ at $c$. Let $c^\p$ be the image of $c$ under $\phi^\p$. Then the isotopy for $t_0-2\le z\le t_0-1-\delta$ is chosen so that the part of $f$ between $b$ and $c$ (which is a vertical band) is mapped to the band joining $b$ and $c^\p$ (which may no longer be a vertical band). This can be done since the shrinkage $\phi^\p$ on $R_0\times \{t_0-2\}$ ($R_0$ is the base of $B_0$) will not cause these bands to intersect each other by the conditions on the positions of the blue squares. This isotopy restores the original boxes $B_1$ and $B_2$ together with all the blue squares.  The blue triangles in $B_0^\p$ remain the same, the vertical blue bands from $z=t_0-1$ to $z=t_0-2$ are no longer vertical (however they remain disjoint from each other and are strictly increasing in the $z$-coordinates). The part of a blue band under $z=t_0-2$ remain vertical after this isotopy is applied.

\medskip
(2) After isotopy (1) the graphs $G_1$ and $G_2$ fit back into the original $B_1$ and $B_2$. The resulting embeddings ``almost" restore $G^g_1$ and $G^g_2$. The ``almost" stems from the exception that the blue triangles in each box that would be part of a grid-like embedding are distorted, however they are identical to $G^g_1$ and $G^g_2$ below the $z$-level where the bases of their blue triangles are. By a slight abuses of notation we call these ``almost" grid-like embeddings still $G_1^g$ and $G_2^g$.  By the definition of grid-like embedding, $G_1^g$ ($G_2^g$) is isotopic to a standard 3D-embedding $G^s_1$ ($G^s_2$) by a VNP-isotopy that is identity outside the box $R_1\times [-1,t_1-\delta]$ ($R_2\times  [-1,t_2-\delta]$ ), see Definition \ref{gridlikedrawing} (v). We will now apply these two isotopies to $G_1^g$ and $G_2^g$. After this, an edge path from the base rectangle $R_0$ to a blue vertex already existed in $G^g_1$ and $G^g_2$ before the reconnection consists of four straight line segments: a single vertical line segment from $R_0$ to $z=t_0-2$, then a line segment (that is in a deformed blue band) from $z=t_0-2$ to $z=t_0-1$, then a vertical line segment from $z=t_0-1$ to the base of a blue triangle in $z=t_0-\delta$, followed a line segment in a blue triangle leading to the blue vertex. Notice that the last two line segments are not changed by the isotopy applied so far.
Note also that these paths do not intersect the connecting rectangle (which is also not affected by the last two isotopies since it is in $z=t_c=t_0-2$). See Figure \ref{deformedband}.

\begin{figure}[htbp]
\begin{center}
\includegraphics[scale=0.6]{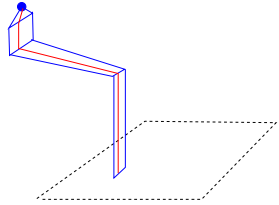}
\caption{An edge path from the base rectangle $R_1$ (or $R_2$) to a blue vertex after the isotopy in (2) is applied to the grid-like embedding $G_0^g$ constructed from $G_1^g$ and $G_2^g$. }
\label{deformedband}
\end{center}
\end{figure}

\medskip
(3) During an edge-cut it is possible that in the contracting process we created temporary loop edges that were deleted in order to form $G_1$ and $G_2$. These loop edges are realized by some ``small hooks" that are attached at the two ends of the neighborhood rectangle on the connecting rectangle ($N^{\p\p}$) below $z=t_0-2$ in the reconnecting process (Step (6) of Section \ref{sec72}). At this stage we must realize these loops in the plane $z=0$ as it is required in a standard 3D-embedding of $G_0$. A top view of two nested such hooks and their relative positions with the other edges involved in the same blue square (on $z=0$) are illustrated in Figure \ref{hookmove}.

\begin{figure}[htbp]
\begin{center}
\includegraphics[scale=0.6]{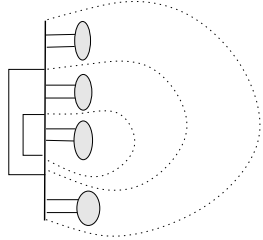}
\caption{A top view of two nested hooks and their relative positions with the other edges involved in the same blue square: the thick line segment represents the segment $L$ used in a blue square, the large outside dashed curve is a conceptual depiction of the deformed $\gamma^\p$, the two small dashed curves represent the actual loops corresponding to the two hooks. }
\label{hookmove}
\end{center}
\end{figure}

\medskip
We accomplish this one hook at a time starting with an innermost hook. We slide a hook along the blue bands down into the plane $z=0$, see Figure \ref{loopisotopy}. Once they are in the plane $z=0$ we fold them by a 90 degree turn into the unused space in the blue square so they look just as shown in Figure \ref{hookmove}. It is clear from the figure that the hook can then be deformed to the dashed curves from under the plane $z=0$ by a VNP-isotopy that is identity below $z=-1$. Of course it needs to remain in the box $B_i$ to which it belongs. See Figure \ref{loopisotopy} for an illustration of this process.

\begin{figure}[htbp]
\begin{center}
\includegraphics[scale=0.6]{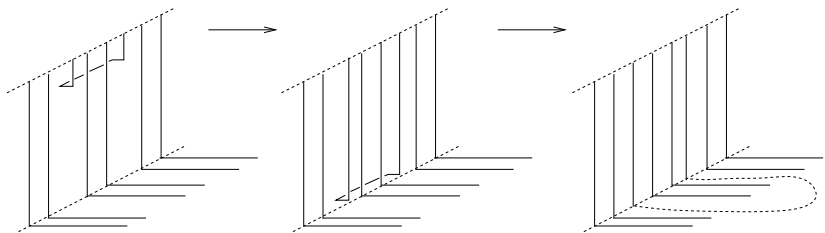}
\caption{How a small hook (loop) is isotoped along a blue band into the plane $z=0$. The edge on the right is dashed because it might have to be stretched out to fit into the plane $z=0$ and can be quite long.}
\label{loopisotopy}
\end{center}
\end{figure}

\medskip
(4) Let us recall that at this point all red vertices are in the plane $z=0$ and all blue vertices are in the plane $z=t_0+2$. Also, all the red-red edges are also in the plane $z=0$ except those going through the connecting strip (created in the re-connecting process). We cannot simply project these into the plane $z=0$ since that will likely to create crossings so we will not be able to recover our original plane graph this way. Since we know that the edges connected to $v_2^b$ in $B_2$ are in the outer face of $G_1$, we will try to shrink the graph $G_1$ first and then move the whole graph along the connecting strip into the interior face $F$ of $G_2$ where it becomes clear that the original graph structure is recovered. The blue bands remain connected to the top of the box $B_0$.  In this way we avoid the creation of unwanted intersections. In the following we describe these steps in more detail.

\medskip
Assume that $w_1$, $w_2$, ..., $w_k$ are the centers of the blue squares associated with the blue vertices in $G_1^g$ other than $v_1^b$. Let $R_s$ be the half of the square $R_{v^b_1}$ in the center of the blue square $S^b_{v_1^b}$ between $L_1$ and $L_{v_1^b}$ as shown in Figure \ref{bandex2}. Choose $k$ points $y_1$, $y_2$, ..., $y_k$ in the small rectangle $R_s$ such that the $y_j$'s do not share the same $x$-coordinates nor $y$-coordinates. This small rectangle can be viewed from the top as the rectangle with dotted line boundary and the letter $G_1^s$ as marked in Figure \ref{bandex2}. By Lemma \ref{redrawGinR}, there exists a plane isotopy $\xi: R_1\times[0,1]\rightarrow R_1$ that is identity outside a small neighborhood of $R_1$ and the identity on $L_{v_1^1}$ (the base of the blue band) that takes $w_j$ to $y_j$ and moves the all points of the embedding $G_2^s$ in the plane $z=0$ into $R_s$. Here we assume that $\xi_0=\xi(x,0)$ is the identity on $R_1$ and $\xi_1=\xi(x,1)$ has moved all the $w_j$'s to the $y_j$'s. This plane isotopy is extended to a VNP-isotopy in the following way: (a) it is the identity outside of the box $R_1\times [-1,t_1-1/2]$;
(b) its action on $R_1\times \{s\}$ for each $0\le s\le t_1-1/2$ is the same as that of $\xi(x,m(1-s/(t_1-1/2)))$ for $m\in [0,1]$ on $R_1\times \{0\}$.
Note that this extension keeps all blue bands disjoint from each other and each edge on a blue band is a path that is non-decreasing in the $z$ direction. Furthermore, the isotopy can be so chosen that the bases of the blue bands are mapped to bands perpendicular to the direction of the extension vector. See Figure \ref{shrink1}.

\begin{figure}[htbp]
\begin{center}
\includegraphics[scale=0.6]{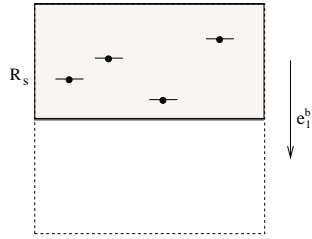}
\caption{How $G_1^g$ looks after it is shrank to fit in $R_s$. The dots are the blue vertices and the bars on them indicate the bases of the blue squares after the isotopy.}
\label{shrink1}
\end{center}
\end{figure}

\medskip
We now define an isotopy that retracts the connecting strip while dragging the rectangle $R_s$ with the edges of $G_1$ along the connecting strip. At first, the rectangle $R_s$ containing $G_1$ can be lifted vertically to the level $z=t_c$.  In doing so we retract the connecting strip at the same time.  All the blue bands that are connecting $R_s$ to the blue vertices at the top of the box $B_0$ are deformed as well to keep the bands nondecreasing in the $z$-coordinates. An illustration of this process is shown in Figures \ref{bandex2} and \ref{bandex4}. This isotopy clearly preserves the neighborhood structures of the vertices.

\begin{figure}[htbp]
\begin{center}
\includegraphics[scale=1.0]{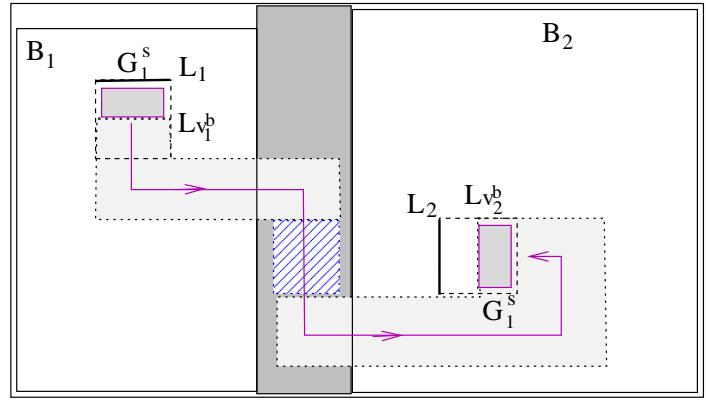}
\caption{Shrinking $G_1$ and moving it to the unused space in the blue square at $v_2^b$.}\label{bandex2}
\end{center}
\end{figure}

\begin{figure}[htbp]
\begin{center}
\includegraphics[scale=1.0]{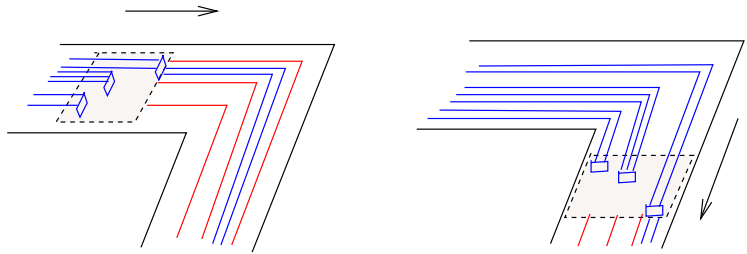}
\caption{Two middle steps of the ``dragging" isotopy: the graph is in the gray area and not explicitly shown to keep the drawing simple. The blue band on the boundary of $R_s$ corresponds to a blue vertex recovered in the reconnection process.}
\label{bandex4}
\end{center}
\end{figure}

\medskip
The retraction continues along the horizontal parts of the extension bands and the connecting rectangle. Whenever the connecting strip turns horizontally we turn the whole rectangle $R_s$ accordingly. Again all the blue bands that are connecting $R_s$ to the blue vertices at the top of the box $B_0$ are dragged along and kept as increasing in their $z$-coordinates as shown in Figure \ref{bandex4}. At the end of the extension band that leads into $B_2$ it is time to drop the small rectangle containing $G_1$ vertically down. By our assumptions on the extension vector at $v^b_2$, there is free space in the small rectangle $R_{v^b_2}$ in the blue square $S^b_{v^b_2}$ to put the box $R_s$ into the plane $z=0$. Once the small rectangle $R_s$ is in the plane $z=0$ the connecting strip has been completely eliminated. We now have the graph $G_0=G(i,j)$ embedded into the plane $z=0$ with the exception of the blue vertices of $G_0$. For each blue vertex there is a small blue square (the blue squares have different sizes) and a blue band that connects the blue square to a blue vertex on the top of $B_0$ along a path that is non-decreasing in the $z$ direction. Since we can deform the blue vertices (with the edges connected to them) into the plane $z=0$ along the deformed bands one by one without any interference and without crossing the boundary of $R_{v^b_2}$, it is clear that the resulting graph bounded within $R_{v^b_2}$ is topologically equivalent to $G_1$ and the resulting graph outside $R_{v^b_2}$ is topologically equivalent to $G_2$. These two graphs are connected along the boundary of $R_{v^b_2}$ following the gluing instruction. Thus the graph obtained after the reconnection is indeed isotopic to $G_0$.

\medskip
(5) Now let $p_1$, $p_2$, ..., $p_j$ be all the intersection points of the edges leading out from the blue vertices of $G^g_0$ with the bases of the blue triangles in $B_0$ (they are all in the plane $z=t_0-\delta$) and let $q_1$, $q_2$, ... $q_j$ be the first intersection points of the corresponding edges with the plane $z=0$. By our construction, $p_i$ and $q_i$ are connected by a path that is non-decreasing in the $z$ direction and these paths do not intersect each other. By Lemma \ref{straighten} and Remark \ref{re22} after it, there exists a VNP-isotopy $\Psi$ such that $\Psi$ is the identity in $z\ge t_0-\delta$ and outside a small neighborhood of $B_0$, and $\Psi$ deforms each path connecting $p_i$ to $q_i$ to a straight line segment and $\Psi$ is a plane isotopy when restricted to $z=0$. The edges in $z=0$ can be further deformed to create spaces for the required blue squares in $z=0$. The result is the desired standard 3D-embedding of $G_0$. This finishes the proof of the first case.

\medskip
Case 2: The case of a vertex-cut subdivision. The situation in this case is slightly simpler because no temporary loop edges are created and the connecting rectangle contains exactly one blue square $S^b_{v_0^b}$. Thus one of the steps in the prior isotopy is no longer needed.
The isotopy is again constructed by a sequence of isotopies using the following steps:

\medskip\noindent
(1) Remove some gridlines to get back to an almost grid-like embedding;

\smallskip\noindent
(2) Shrink $B_1^\p$ and $B_2^\p$ (with the graphs contained in them) back to $B_1$, $B_2$ and deform the graphs contained in them to 3D standard embeddings;

\smallskip\noindent
(3) Shrink $G_1$ and drag the shrank $G_1$ along the extension bands and the connecting rectangle and drop it into $G_2$ on $z=0$.

\smallskip\noindent
(4) Straighten out the blue squares and make the blue-bands vertical as required by a standard 3D-embedding of $G_0$ using Lemma \ref{straighten}.

\medskip
Note that it is essential that $G_1$ is used in step (4). For example, suppose we are in case (ii) of Lemma \ref{diskcomponent}, that is $G_1$ consists of a union of several maximal disk-components of a single inseparable disk-component. Then it is possible to drag $G_1$ into $G_2$ in step (3) above, but it might not be possible to drag $G_2$ into $G_1$ without getting hung up on some blue band that is connected to the top of the box $B_0$.

\medskip
Some more argument is needed for the isotopy in Step (4). In the case for the edge-cut, the shrunk graph $G_1$ is positioned in the unused rectangle in $R_{v^b_2}$ so placing $G_1$ into that space apparently will not cause of edge intersections. However, in this case the shrunk graph $G_1$ is positioned into the {\it used} half of $R_{v^b_2}$ in $G_2$. Since $G_1$ contains some (or all) maximal disk components of some inseparable disk component $Q$ of $G_0$, the edges from $G_1$ connected to the blue band $B^b_{v^b_1}$ fit precisely between the edges of the $Q\setminus G_1$ in $G_2$ already connecting to $L_{v^b_2}$. If one combines Figures \ref{principalvertexcut} and \ref{bandex4} (but without the red-red edges shown), then this becomes clearer.

\medskip
The details of the isotopies are otherwise similar to the edge-cut argument and are thus left to the reader. This finishes the proof and also concludes this section.
\end{proof}

\medskip
\section{Embedding Length Analysis}\label{sec9}

\medskip
One way to estimate the embedding length is to estimate the volume of the box $B_0$ in this construction. By our design of the construction, the height of $B_0$ is only two units larger than the maximum height of the two boxes $B_1$ and $B_2$. We can thus estimate the volume of $B_0$ by estimating the area of the base rectangle of $B_0$. One way to do this is to count the number of horizontal gridlines in the $x$- and $y$-direction we must create to accommodate a newly created red-red edge. A gridline is consumed either by a segment of the red-red edge in the $xy$-plane using that gridline or by the red-red edge moving through a vertical segment parallel to the $z$-axis. More precisely we are interested in the following. Let $e$ be a newly created red-red edge and consider a vertical projection $p(e)$ of $e$ into the $xy$-plane. We count the 90 degree turns of $p(e)$ and add to it the number of segments parallel to the $z$-axis to obtain an estimate on the number of gridlines that need to be available in the base rectangle $R_0$.

\begin{definition}
\label{turningnumber}
Let $(G^g, B)$ be a grid-like embedding of a BRT-graph $G$ contained in a rectangular box $B$. Let $E_{br}(G)$ be the set of red-blue edges in $G$. For any edge $e\in E(G)$ ($e$ can be a red-red or a blue-red edge) denote with $tr(e)$ the number of different horizontal gridlines (in the $x$- and $y$-direction) needed and call this the {\em turning number} of the edge $e$. We define the {\em turning number} $T(G^g,B)$ of $(G^g, B)$ as follows:
$$
T(G^g,B)=\max_{e\in E_{br}(G)}\{tr(e)\}.
$$
\end{definition}

Note that in the above definition we ignore the turns which are not a 90 degree angle that may occur when an edge enters a red rectangle or a blue triangle. Counting the 90 degree turns in $p(e)$ over counts the gridlines needed for the horizontal segments in the red-red edge, since several segments of $p(e)$ may end up on the same gridline. Note also that $p(e)$ may have 180 degree turns that arise when two consecutive turns in the $xz$ or $yz$ plane occur, as is the case for a hook. Those are not counted, since the vertical segments are counted separately. We are now able to show the following Theorem.

\begin{theorem}
\label{sumextconstr}
Let  $G_1$ and $G_2$ be the BRT-graphs obtained from the BRT-graph $G_0$ by a single subdivision (either a vertex-cut or an edge-cut). Let $G^g_0$ be the grid-like embedding of $G_0$ in a rectangular box $B_0=R_0\times [0,t_0]$, obtained by reconnecting the grid-like embeddings $G_1^g$, $G_2^g$  of  $G_1$, $G_2$ as described in Section \ref{sec72}.  If $T(G^g_1,B_1)=n_1$ and $T(G^g_2,B_2)=n_2$, then the grid-like embedding $G^g_0$ of $G_0$ has the following properties:

\medskip\noindent
(i) For a newly constructed red-red edge $e$ of $G_0$ $tr(e)\le 2\max \lbrace n_1,n_2 \rbrace  + 8$.

\smallskip\noindent
(ii) $T(G^g_0,B_0)\le \max \lbrace n_1,n_2 \rbrace + 17$.
\end{theorem}

\begin{proof}
We first prove (i). Note that a newly constructed red-red edge can only appear if $G_1$ and $G_2$ are created by an edge-cut in $G_0$. The path a newly created red-red edge $e$ travels is divided into three parts, where each part is a union of one or more line segments each of which is parallel to one of the coordinate axis. There are two cases to consider: the curve $\gamma^\p$ (the push-off of $\gamma$ as defined in Section  \ref{sec3}) intersects $e$ once or the curve $\gamma^\p$ intersects $e$ twice.

\medskip
In the first case when traveling along the edge $e$ we encounter the following parts: The first part starts at a red vertex, say in $G_1^g$, moves to the blue square $S^b_{v^b_1}$, and continues through the blue square to the base of the blue band $B_{v^b_1}$ and then up along the blue band. At this point we know that this first part of $e$ has at most turning number $n_1$. The middle part of $e$ begins with a 90 degree turn in the $xz$- or $yz$-plane, from the vertical direction on the blue band $B_{v^b_1}$ to the extension band, moves through the connecting strip and ends with a 90 degree turn (again in the $xz$- or $yz$-plane) from the horizontal extension band to the vertical direction on the blue band $B_{v^b_2}$. It takes at most 2 turns to reach the connecting rectangle from each box $B_i$ and two turns to turn onto and off the connecting rectangle, giving us a turning number of at most 6 for the middle part. The two 90 degree turns off and onto the vertical blue band $B_{v^b_i}$ are not counted, since the need for this gridline is already accounted for when the blue band is created. Thus $tr(e)\le n_1+n_2+6$ in $G_0^g$.

\medskip
If we assume that an edge $e= w_1w_2$ is intersected by $\gamma^\p$ twice then both vertices $w_1$ and $w_2$ are in the same graph $G_i$.  The middle part of $e$ that does not contain any of the two vertices, is contracted into a loop edge and then deleted. In the reconstruction process a small hook is constructed whose projection looks like the letter $H$. The horizontal bar of the $H$ is the edge segment of the hook in the $xy$-plane below the connecting rectangle which is not connected to a vertical segment and one gridline needs to be added for it. One pair of ends of the two parallel lines in the $H$ is connected to $w_1$ and $w_2$, the other pair of ends is connected to the vertical segments of the constructed hook. Since both ends are on the same line (in the $x$- or $y$-direction) only one gridline needs to be added even though there are two vertical segments.  We account for the fact that both $w_1$ and $w_2$ are in the same graph by using a maximum. In this case  the turning number satisfies $tr(e)\le 2 \max \lbrace n_1,n_2 \rbrace + 8$ turns in $G_0^g$.

\medskip
Next we prove (ii) by proving that $tr(e)$ for a red-blue edge $e$ is $ \le \max \lbrace n_1,n_2 \rbrace + 16$. We consider the case of an edge-cut first. Let $e$ be a red-blue edge and assume that the red vertex $w$ of $e$ is contained in $G^g_1$ (the case of $G^g_2$ is accounted for by using a maximum as before) and that $x$ is the blue vertex of $e$. We note that the curve $\gamma^\p$ can intersect $e$ at most once. The result of the edge-cut is a red-blue edge in $G^g_1$ from $w$ to the blue vertex $v^b_1$ and potentially a deleted loop edge in $G_2$ that results in a hook in $G^g_0$. No loop (and thus no hook) is created if the blue vertex $x$ is on $\gamma^\p$, or if $x$ is in $G_2$ and is not connected to any edge other than $e$ which is cut by $\gamma^\p$.

\medskip
First let us consider the case when there is no hook in $G^g_0$. In this case, $e$ consists  of two parts in $G^g_1$. As before, the first part starts at the red vertex $w$ in $G_1^g$, moves to the blue square $S^b_{v^b_1}$, and continues through the blue square to the base of the blue band $B_{v^b_1}$ and then up along the blue band. This part contains at most $n_1$ turns.
The second part of $e$ begins with a 90 degree from the blue band $B_{v^b_1}$ to the extension band, moves through the extension band, onto the middle section of the connecting rectangle with at most 3 turns, then moves up to the new blue square with at most 4 more turns, once in the blue square, it may need up to another 6 turns to be re-routed to enter the middle square of the blue square in the correct order (to achieve the desired linear order), then finally move up to the base of the new blue triangle with one more turn. This part of the path will add at most 15 right angle turns for the second part (and 14 total) in this case.

\medskip
In case there is a hook,  more turns are needed for the second part. As before, the edge turns from the blue band onto the extension band, moves through the extension band, across the connecting rectangle, runs through the small hook, returns to the middle section of the connecting rectangle and then move up a level to a new blue square and finally to the base of a blue triangle. ends in one of the blue vertices $w_{c_i}$ of $G^g_0$ on the connecting rectangle. The extra number of gridlines needed for the loop is 2, as already established during our earlier discussion for the red-red edges. Thus $tr(e)\le \max \lbrace n_1,n_2 \rbrace + 17$ for an edge-cut.

\medskip
We now consider the case of a vertex-cut. Here the situation is easier since the red-blue edge $e$ is not cut into parts at all. It already exists in one of the two graphs $G^g_i$, $i=1$, $2$ and it is simply extended from the blue band $B_{v^b_i}$ to the single blue vertex $v^b$ on the connecting rectangle. It is easy to see that this extension does not exceed $17$ turns and thus the results from the edge-cut case suffices as an upper bound for $tr(e)$.
\end{proof}

\medskip
We are now ready to prove the main theorem of this paper as stated below.

\begin{theorem}
\label{mainT}
Let $G$ be a 4-regular plane graph with $n$ vertices. Then there exists a realization of $G$ on the cubic lattice $\Z^3$ which is contained in a rectangular box whose volume is bounded above by $ O(n\ln^5n)$. Consequently, the ropelength of any knot or link $K$ is at most of the order $O(Cr(K)\ln^5(Cr(K))$, where $Cr(K)$ is the minimum crossing number of $K$.
\end{theorem}

\medskip
To prove the theorem, we will use the non-linear recurrence analysis approach as described in \cite{Lei}. The result of Theorem \ref{sumextconstr} and the results from the earlier sections will be needed when we apply this analysis.

\medskip
For any rectangle $R$ under discussion we assume that its
length $l$ is greater or equal to its width $w$. The {\it aspect
ratio} $\sigma$ of the rectangle $R$ is defined by $\sigma = w/l$. From
this it follows that $l=\sqrt{A/\sigma}$ and $w=\sqrt{\sigma A}$
where $A$ is the area of the rectangle $R$. The aspect ratio is
important to us because throughout the algorithm we want to operate
with boxes whose base rectangles are not too skinny, that is the
rectangles have an aspect ratio that is bounded away from zero by a positive constant.

\medskip
In the previously described divide-and-conquer algorithm we divide a
given BRT-graph at each step into subgraphs (also
BRT-graphs) each of which has a size of at least 1/6
of the previous graph, where the size is measured by the number of red
vertices in a graph.  This allows us to operate with rectangles whose
aspect ratios are at least 1/6. Before getting started on the details
we need a preliminary lemma that asserts that we can divide
rectangles while preserving the minimal aspect ratio of 1/6.

\begin{lemma}
\label{aspectratio}
Let $R$ be a rectangle with length $l$ and width $w$ and aspect ratio
$\sigma =w/l\ge 1/6$. Let $1/6\le\alpha\le 5/6$ be a real number and
divide the rectangle $R$ by a line that is parallel to its width into
two rectangles $R_1$ and $R_2$ of areas $\alpha w l$, $(1-\alpha) w
l$ and aspect ratios $\sigma_1$, $\sigma_2$ respectively, then $\sigma_i\ge 1/6$ for each $i$.
\end{lemma}

\begin{proof}
It suffices to show this for one of the aspect ratios, say $\sigma_1$ of the rectangle $R_1$. We
can assume that the area of $R_1=\alpha w l$ with dimensions $\alpha l$ and $w$. There are two cases to
consider,  either $w\ge \alpha l$ or $w< \alpha l$. If $w\ge \alpha l$
then $\sigma_1=\alpha/\sigma$ and $1/6\le(1/6)/\sigma\le \alpha/\sigma=\sigma_1$.
If $w< \alpha l$
then $\sigma_1=\sigma/\alpha$ and $1/6<(1/6)/\alpha\le \sigma/\alpha=\sigma_1$.
\end{proof}

From the previous sections we know that a BRT-graph $G$ with $n$ red vertices has a balanced recursive subdivision sequence whose depth is bounded above by $c_r \ln (n)$ for some constant $c_r > 0$, see Theorem \ref{Trecexist}. At the end of the subdivision process we have terminal BRT-graphs $G(i,j)$ that satisfy $W_s(G(i,j)) < W_0$. That is, each such terminal BRT-graph has less than $W_0$ red vertices. By Lemma \ref{vbbound}, the number of blue vertices in $G(i,j)$ is bounded above by $3W_0$ (since the $g$ used in Lemma \ref{vbbound} is equal to 12). Thus all graphs $G(i,j)$ have at most $n_0=4 W_0$ vertices (blue and red combined). Since there are only finitely many plane BRT-graphs with at most $W_0$ red vertices and whose maximum vertex degree is 12 or less, the following lemma holds for a grid-like embedding $G^g(i,j)$ of a terminal BRT-graph $G(i,j)$.

\medskip
\begin{lemma}\label{N0lemma}
There exists an integer $N_0>0$ such that if $G(i,j)$ is a terminal BRT-graph (so it has less than $W_0$ red vertices) and $R$ is any rectangle of area $\ge N_0$ with an aspect ratio $\sigma\le 1/6$, then $R$ is large enough such that the rectangular box $B(i,j)$ of height 1 with $R$ as its base rectangle can be used to hold a grid-like embedding $G^g(i,j)$ of $G(i,j)$. Moreover there exists a constant $m>0$ that only depends on $W_0$ such that for any such grid-like embedding $G^g(i,j)$ in the box $B(i,j)$, the turning number $T(G^g(i,j),B(i,j))\le m$.
\end{lemma}

\medskip
Let us now examine what happens in the reconstruction process. Let $i_0$ be the depth of the subdivision sequence. Recall that when two grid-like embeddings are reconnected, the turning number of a red-blue edge that remains a red-blue edge increases by at most $c_t$ for some constant $c_t$ ($c_t$ can be chosen to be $17$ in fact, see Theorem \ref{sumextconstr} (ii)).

Therefore if a grid-like graph $G^g(i_0-1,j)$ is reconstructed from the grid-like embeddings $G^g(i_0,2j-1)$ and $G^g(i_0,2j)$ in a box $B(i_0-1,j)$, the height of $B(i_0-1,j)$ is $1+2=3$ and $T(G^g(i_0-1,j),B(i_0-1,j))\le m+c_t$. Inductively, we can see that if a grid-like graph $G^g(i,j)$ is reconstructed from the grid-like embeddings $G^g(i+1,2j-1)$ and $G^g(i+1,2j)$ in a box $B(i,j)$, then the height of $B(i,j)$ is $1+2(i_0-i)$ and $T(G^g(i,j),B(i,j))\le m+c_t(i_0-i)$ for any $0\le i\le i_0-1$. By Remark \ref{depthbound}, we can then claim the following.

\medskip\noindent
(i) The height of the box $B(i,j)$ is bounded above by $1+2c_r \ln(W_s(G(i,j))$;

\smallskip\noindent
(ii) $T(G^g(i,j),B(i,j))$ is bounded above by $m+c_t\ln(W_s(G(i,j)))$.

\medskip
If a red-red edge $e$ is created in $G^g(i,j)$ by connecting two red-blue edges (one in $G^g(i_0,2j-1)$ and one in $G^g(i_0,2j)$), then the turning number of the resulting edge $e$ is the sum of the two turning numbers of the two red-blue edges in $G^g(i_0,2j-1)$ and $G^g(i_0,2j)$, as well as an addition of at most 6, as shown in Theorem \ref{sumextconstr} (i). By the inequality given in (ii) above, the turning number of the edge $e$ is then bounded above by $2m+2c_t(i_0-i-1)+6$. However the turning number of $e$ no longer changes in the subsequent reconnection process. We can summarize this as follows:

\medskip
(iii) Let $G(i,j)$ be a BRT-graph that is subdivided into $G(i+1,2j-1)$ and $G(i+1,2j)$. Then for any newly created red-red edge $e$ in the construction process of $G^g(i,j)$ we have
$$
tr(e)\le 2m+6+2c_t(c_r\ln(W_s(G(i,j)))-1).
$$

\medskip
In the reconstruction of $G^g(i,j)$ from $G^g(i+1, 2j-1)$ and $G^g(i+1, 2j)$, for each new red-red edge $e$ created, we need to add up to $tr(e)$ many gridlines to place $e$ on the lattice. The number of new red-red edges created at that stage is the same as the number of red-red edges cut by $\gamma^\p$. Since the number of edges cut by $\gamma^\p$ is proportional to the square root of the total number of the vertices in the BRT-graph $G(i,j)$, by Lemma \ref{vbbound}, the number of edges cut is proportional to the square root of the red vertices of $G(i,j)$. Thus the total number of gridlines which may have to be added to create the rectangle $R(i,j)$ from the rectangles  $R(i+1,2j-1)$ and $R(i+1,2j)$ is at most

\begin{equation}\label{bound_const}
c \cdot(\ln(W_s(G(i,j)))\sqrt{W_s(G(i,j)})
\end{equation}

for some constant $c>0$ that only depends on the construction algorithm (by the inequalities we obtained from (i)--(iii) above).

\medskip
\begin{definition}
Let $A(n)$ be the minimum over all positive
integers $p$ with the property that if the area of a rectangle
$R$ with aspect ratio $\sigma_R\ge 1/6$ is equal
to or greater than $p$, then any BRT-graph $G$ with $n$ red
vertices whose degrees are at most $12$ has a grid like embedding in a
box $B=R\times [0, c_h]$ where $c_h$ is a positive integer and $c_h\le 2\lceil\frac{\ln(n)}{\ln(6/5)}\rceil+1$ .
\end{definition}

\medskip
In the following theorems $\alpha$ or $\alpha_i$ is always a
number in [1/6, 5/6]. However only those values in [1/6, 5/6] whose product with $n$ is an integer make sense where $n$ is the number of red vertices of the graph that is subdivided. The reader should keep this in mind in the statement and proof of Theorem \ref{upperboundAn} below.

\begin{theorem}
\label{upperboundAn}
There exists a constant integer $N_0$ such that the function $A^\p(n)$ defined by
$$
A^\p(n)=\left\{
\begin{array}{ll}
N_0 & {\rm if}\ n\le W_0,\\
\max_{1/6\le\alpha\le 5/6} (\sqrt{A^\p(\alpha
n)+A^\p((1-\alpha)n)}+c\sqrt{6n}\ln(n))^2 & {\rm if}\ n>W_0
\end{array}
\right.
$$
bounds $A(n)$ above, where $c$ is the same constant used in Equation \ref{bound_const}.
\end{theorem}

\begin{proof}
We will prove this theorem inductively. We choose $N_0$ as the constant guaranteed by Lemma \ref{N0lemma}. By definition, we have $A(k)\le A^\p(k)=N_0$ if $k\le W_0$. Now assume that $A(k)\le A^\p(k)$ is true for any $k\le n-1$ for $n\ge W_0+1$, we need to prove that $A(n)\le A^\p(n)$.

\medskip
Let $\sigma\le 1/6$ be any given aspect ratio and $R$ a rectangle with area $A^\p(n)$ and an aspect ratio $\sigma$. Let $G$ be a BRT-graph with $n$ red vertices. We know that there exists a subdivision of $G$ that divides $G$ into two
BRT-graphs $G_1$ and $G_2$ with $\alpha n$ and $(1-\alpha)n$
red vertices respectively for some $\alpha\in [1/6, 5/6]$.
Let $R_0$ be a rectangle with aspect ratio $\sigma$ and area
$A^\p(\alpha n)+A^\p((1-\alpha)n)$.
By Lemma \ref{aspectratio}, $R_0$ can be divided into two rectangles $R_1$ and $R_2$ whose aspect ratios are in $[1/6,1]$ and whose areas are $A^\p(\alpha n)$ and $A^\p((1-\alpha)n)$ respectively. Thus by our induction hypothesis, there exist grid-like embeddings $G^g_1$, $G^g_2$ of $G_1$, $G_2$ respectively in the rectangular boxes $B_1$, $B_2$ whose bases are $R_1$, $R_2$ and whose heights are at most $2\lceil\frac{\ln(\alpha n)}{\ln(6/5)}\rceil+1$, $2\lceil\frac{\ln((1-\alpha)n)}{\ln(6/5)}\rceil+1$ respectively. Suppose that $R_0$ is placed in such a way that its longer side is parallel to the $x$-axis as shown in Figure \ref{rectangle}. We will now extend $R_0$ to a new rectangle $R^\p$ by adding a strip of width
$c \sqrt{n} \ln (n)$ to its top and adding a strip of width $(c/\sigma) \sqrt{n} \ln (n)$ to its right. The strip added to the right of $R_0$ is thicker than the strip added to its top since $c/\sigma>c$. Observe that $R^\p$ also has aspect ratio $\sigma$.

\begin{figure}[htbp]
\begin{center}
\includegraphics[scale=0.5]{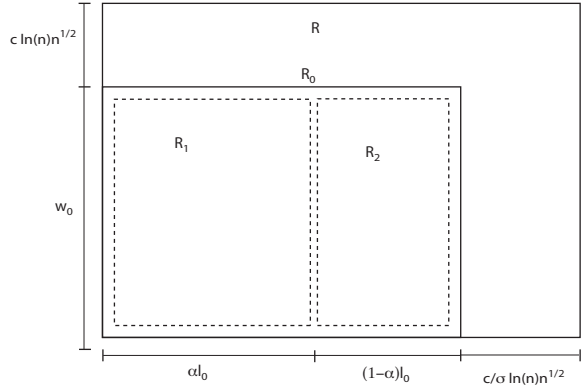}
\caption{The rectangles $R_0$, $R_1$, $R_2$ and $R^\p$. The rectangles $R_1$ and $R_2$ are the base rectangles of the boxes $B_1$ and $B_2$. }\label{rectangle}
\end{center}
\end{figure}

\medskip
We know that we can obtain a grid-like embedding of $G$ from $G_1^g$ (which is bounded in $B_1$ with $R_1$ as its base) and $G_2^g$ (which is bounded in $B_2$ with $R_2$ as its base). At most $c\sqrt{n}\ln (n)$ new horizontal and $c \sqrt{n} \ln (n)$ new vertical grid lines are needed to accommodate the red-red edges (which need to be put on the lattice). This means that the base of the new rectangular box $B$ (that houses this new grid-like embedding of $G$) can be fit into $R^\p$. Furthermore, the height of $B$ is 2 plus the larger of the heights of $B_1$ and $B_2$. Say the height of $B_1$ is larger. Then the height of $B$ is at most $2\lceil\frac{\ln(\alpha n)}{\ln(6/5)}\rceil+3$. Since $\alpha\le 5/6$, $\ln(\alpha n)\le \ln n+\ln (5/6)=\ln n-\ln(6/5)$. Thus the height of $B$ is at most $2\lceil\frac{\ln(\alpha n)}{\ln(6/5)}\rceil+3\le 2\lceil\frac{\ln(n)}{\ln(6/5)}-1\rceil+3=2\lceil\frac{\ln( n)}{\ln(6/5)}\rceil+1$. In other words, $G$ has a grid-like embedding in the box $B^\p$ with $R^\p$ as its base and with a height at most $2\lceil\frac{\ln( n)}{\ln(6/5)}\rceil+1$.
Let $w_0$ and $\ell_0$ be the width and length of $R_0$ respectively.
Then $\ell_0=w_0/\sigma$ and $A^\p(\alpha n) +A^\p((1-\alpha) n)=\ell_0w_0=w_0^2/\sigma$, so $w_0=\sqrt{\sigma(A^\p(\alpha n) +A^\p((1-\alpha) n))}$. Thus the area of $R^\p$ is
\begin{eqnarray*}
&& (w_0/\sigma + (c/\sigma)\sqrt{n}\ln(n))(w_0+c\sqrt{n}\ln(n))\\
&=& (w_0+c\sqrt{n}\ln(n))^2/\sigma\\
&=& (\sqrt{\sigma(A^\p(\alpha n) +A^\p((1-\alpha) n))}+c\sqrt{n}\ln(n))^2/\sigma\\
&=& (\sqrt{A^\p(\alpha n) +A^\p((1-\alpha) n)}+c\sqrt{n/\sigma}\ln(n))^2.
\end{eqnarray*}

However since $1/6\le\sigma$, $n/\sigma\le 6n$, the area of $R^\p$ is bounded above by
\begin{eqnarray*}
&&(\sqrt{A^\p(\alpha n) +A^\p((1-\alpha) n)}+c\sqrt{6n}\ln(n))^2\\
&\le & \max_{1/6\le\alpha\le 5/6} (\sqrt{A^\p(\alpha
n)+A^\p((1-\alpha)n)}+c\sqrt{6n}\ln(n))^2= A^\p(n).
\end{eqnarray*}

Since $R$ and $R^\p$ have the same aspect ratio and $R$ has a larger area, $R^\p$ can be fit into $R$. Therefore, $G$ has a grid-like embedding in a box with base $R$ and height at most $2\lceil\frac{\ln( n)}{\ln(6/5)}\rceil+1$. This proves that $A(n)\le A^\p(n)$.
\end{proof}

\medskip
The next theorem gives the function $A^\p(n)$ (hence $A(n)$) an explicit bound.

\begin{theorem}
There exists a constant $d>0$ such that
$A^\p(n)\le d n (\ln (n))^4.$ It follows that $A(n)\le d n (\ln (n))^4$ as well.
\end{theorem}

\begin{proof}
Following the proof given in \cite{Lei} we define a function $B(n)$
as follows: For $n \le W_0$, $B(n)=\sqrt{N_0}$ and for $n > W_0$
\begin{eqnarray}
\label{eqn1}
B(n)= \max_{1/6\le\alpha\le 5/6} (B(\alpha n)+c\sqrt{6}\ln(n)).
\end{eqnarray}

We show that $A^\p(n)\le n (B(n))^2$ for all $n\ge W_0$ by induction.
Clearly this is true for $n=W_0$. Assume that $A^\p(k)\le n (B(k))^2$ is true for all values of $k$ such that $W_0\le k<n$. For $k=n$ we have:

\begin{eqnarray*}
A^\p(n)&=&\max_{1/6\le\alpha\le 5/6} (\sqrt{A^\p(\alpha n)+A^\p((1-\alpha)n)}+c\sqrt{6n}\ln(n))^2
\\
&\le&\max_{1/6\le\alpha\le 5/6} (\sqrt{\alpha n (B(\alpha
n))^2+(1-\alpha)n(B((1-\alpha) n))^2}+c\sqrt{6n}\ln(n))^2
\\
&\le&\max_{1/6\le\alpha\le 5/6} (\sqrt{ n (B(\alpha n))^2}+c\sqrt{6n}\ln(n))^2
\\
&\le&\max_{1/6\le\alpha\le 5/6} n(B(\alpha n)+c\sqrt{6}\ln(n))^2
=n (B(n))^2.
\end{eqnarray*}

The third line in the above inequalities can be explained as follows:
If $B(\alpha n)\ge B((1-\alpha)n)$ for the value of $\alpha$
realizing the maximum then this is obvious. If  $B(\alpha n)<
B((1-\alpha)n)$ then the same result follows where $\alpha$ is
replaced by $(1-\alpha)$. Then a change of variable of
$(1-\alpha)$ for $\alpha$ produces the same result.

\medskip
Now $B(n)$ needs to be estimated. Using Equation (\ref{eqn1})
repeatedly for different values of $\alpha$ results in:
$$
B(n)= \sqrt{6}c(\ln(n)+\ln(\alpha_1n)+\ln(\alpha_1\alpha_2n)+\ldots +
\ln(\alpha_1\ldots\alpha_s n)+N_0,
$$
where $s$ is the depth of the recursion and each value $\alpha_i$ is
the value of $\alpha$ that realizes the maximum at each stage of the
recursion. Since all $\alpha_i\le 1$ we have for some constant $c^\p$
$$
B(n)\le \sqrt{6}c\cdot s(\ln(n))+N_0\le c^\p (\ln(n))^2.
$$
From this it follows that there exists a constant $d>0$ such that
$A^\p(n)\le n (B(n))^2\le d \cdot n (\ln(n))^4$.
\end{proof}

\begin{corollary}
\label{graphlatticeembed}
Let $G$ be a 4-regular plane graph with $n$ vertices. Then there exists a realization of $G$ on the cubic lattice $\Z^3$ which is contained in a rectangular box whose volume is bounded above by $ O(n\ln^5n)$.
\end{corollary}

\begin{proof}
First, we change $G$ to a BRT-graph $G^\p$ by triangulation. Next we construct a grid-like embedding of $G^\p$ in the lattice in a box $B=R\times I$ by the previously described algorithm. The number of vertices in $G$ is equal to the number of red vertices in $G^\p$. The algorithm described in the paper generates a rectangular box with a height of at most $c_h\ln n$ for some constant $c_h>0$. From Theorem \ref{upperboundAn} we know that the area of the rectangle $R$ is of the order $O(n\ln^4n)$ and thus the volume of the box containing the grid-like embedding of $G$ is bounded above by $O(n\ln^5n)$.

\medskip
The lattice embedding described in the prior section assumes that edge segments in the red rectangles are not on lattice. In order to obtain a complete lattice embedding of $G$, the edges which are added to $G$ in the triangulation to obtain the BRT-graph $G^\p$ are removed from the embedding. Now each red rectangle only intersects 4 edges and these edges can be re-arranged in the rectangle such that they connect to the red vertex using only lattice connections.  Thus the lattice embedding of $G$ fits into the same rectangular box as the grid-like embedding of $G$.
\end{proof}

\begin{figure}[htbp]
\begin{center}
\includegraphics[scale=0.5]{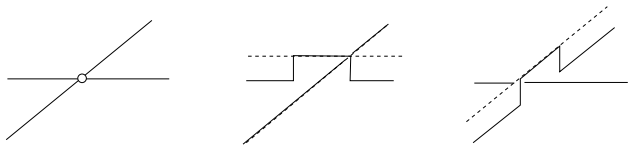}
\caption{Reroute the edges at a vertex to recover the desired under/over crossing. }\label{crossingrecover}
\end{center}
\end{figure}

Let $D$ be a knot diagram of $\K$ with $n$ crossings. $D$ can be thought of as a 4-regular plane graph and by Corollary \ref{graphlatticeembed} there exists a lattice embedding of the 4-regular plane graph $D$ in a rectangular box with a volume of order $O(n\ln^5n)$. In order to change this lattice embedding of the graph $D$ to the embedding of the knot $\K$, we need to recover the crossings of the original knot diagram from the vertices of the embedded graph. In order to accomplish this, the edges at a vertex $v$ are locally modified as shown in Figure \ref{crossingrecover}. One unit away from the vertex $v$ the edge is rerouted through the plane $z=0.5$ to recover the desired crossing. The space between $z=0$ and $z=0.5$ and the space between $z=0.5$ and $z=1$ are then stretched to one unit thick each (so that $z=0.5$ becomes the lattice plane $z=1$ afterwards). The new rectangular box containing the recovered knot is identical to the previous one except that its height has increased by one. The volume of the surrounding box is still bounded above by $O(n\ln^5n)$.

\medskip
The length of the lattice embedding of a diagram $D$ is bounded above by the volume of the box surrounding the lattice embedding of $D$. Furthermore, the ropelength of a knot $\K$ is bounded above by twice the length of the lattice embedding. The result of our main theorem (Theorem \ref{mainT}) then follows.

\section{Some open questions}\label{sec91}
\medskip
Theorem \ref{mainT}  answers the following question raised in \cite{DE2,DE3} negatively.

\begin{step}
{\bf Question}. For any $1<p\le 3/2$, does there exist
a family of infinitely many knots (links) such that  $L(K)\ge O(Cr(K))^p)$ for knots in this family?
\end{step}

\medskip
However the following question also raised in \cite{DE2,DE3} remains open:

\begin{step}
{\bf Question}. Is it true that $\sup\{{L(K)\over Cr(K)}\}=\infty$
(where the supremum is taken over all knots and links)?
\end{step}

\medskip
Equivalently we can ask the following:
\begin{step}
{\bf Question}. Does there exist a $p>0$ and a constant $a>0$ such that there exists an infinite family of knots and links such that for any member
$K$ in the family, $L(K)\ge a\cdot (Cr(K))\cdot \ln(Cr(K))^p$?
\end{step}

Since we know that there exist knots $\K$ with ropelength of order $O(Cr(\K))$, the only possible further improvement on Theorem \ref{mainT} is to reduce the power on the $\ln(Cr(K))$ term. In \cite{DEZ}, a numerical study suggests the possibility of a ropelength upper bound of the form $O(Cr(K)\ln^2(Cr(K)))$. We thus end this paper with the following question:

\medskip
\begin{step}
{\bf Question}. Can we improve the ropelength upper bound to $L(K)\ge a\cdot (Cr(K))\cdot \ln(Cr(K))^2$?
\end{step}

\medskip
\bibliographystyle{amsalpha}

\end{document}